\documentclass[11pt]{amsart}
\usepackage[margin=1in]{geometry}

\usepackage{amssymb}
\usepackage{amsthm}
\usepackage{amsmath}
\usepackage{mathrsfs}
\usepackage{amsbsy}
\usepackage{bm}
\usepackage{mathtools}
\usepackage{graphicx}
\usepackage{tikz}
\usepackage{array}
\usepackage{enumerate}
\usepackage{xcolor}
\usepackage{comment}
\usepackage{booktabs}
\usepackage[expansion=false]{microtype}
\usepackage[noadjust]{cite}
\usepackage{hyperref}
\usepackage{caption}
\usepackage[noabbrev,capitalise,nameinlink]{cleveref}

\definecolor{lavender}{rgb}{0.4,0,1}
\definecolor{Teal}{rgb}{0,0.60,0.60}
\definecolor{LightBlue}{rgb}{0.392,0.392,1}
\definecolor{DarkPink}{rgb}{0.85,0,0.38}
\hypersetup{colorlinks=true,citecolor=lavender,linkcolor=lavender,urlcolor=lavender}
\allowdisplaybreaks

\newtheorem{theorem}{Theorem}[section]
\newtheorem{proposition}[theorem]{Proposition}
\newtheorem{corollary}[theorem]{Corollary}
\newtheorem{conjecture}[theorem]{Conjecture}
\newtheorem{question}[theorem]{Question}

\newtheorem{lemma}[theorem]{Lemma}

\theoremstyle{definition}

\newtheorem{example}[theorem]{Example}
\crefname{conjecture}{Conjecture}{Conjectures}
\crefname{question}{Question}{Questions}

\newcommand{\dfn}[1]{\textcolor{Teal}{\emph{#1}}}

\newcommand{\SSS}{S}
\newcommand{\ZZ}{\mathbb Z}
\newcommand{\RR}{\mathbb R}
\newcommand{\CC}{\mathbb C}
\newcommand{\NN}{\mathbb Z_{\geq0}}
\newcommand{\one}{\mathbf 1}
\newcommand{\xx}{\mathbf x}
\newcommand{\yy}{\mathbf y}
\newcommand{\zz}{\mathbf z}
\newcommand{\ww}{\mathbf w}
\renewcommand{\aa}{\mathbf a}
\newcommand{\qq}{\mathbf q}
\newcommand{\cF}{\mathcal F}
\newcommand{\cB}{\mathcal B}

\newcommand{\cV}{\mathcal V}
\newcommand{\cN}{\mathcal N}

\newcommand{\cD}{\mathscr D}
\newcommand{\cA}{\mathcal A}
\newcommand{\Hess}{\operatorname{Hess}}

\newcommand{\conv}{\operatorname{conv}}
\newcommand{\rk}{\operatorname{rk}}
\newcommand{\Des}{\mathrm{Des}}
\newcommand{\des}{\mathrm{des}}
\newcommand{\pk}{\operatorname{peak}}
\newcommand{\DT}{\operatorname{DT}}
\newcommand{\AT}{\operatorname{AT}}
\newcommand{\NC}{\operatorname{NC}}
\newcommand{\Fer}{\operatorname{Fer}}
\newcommand{\indeg}{\operatorname{indeg}}
\newcommand{\Sym}{\operatorname{Sym}}
\newcommand{\wt}{\operatorname{wt}}
\newcommand{\black}{\operatorname{black}}
\newcommand{\white}{\operatorname{white}}
\newcommand{\dd}{\,\mathrm d}
\newcommand{\halfplane}{\mathbb H}
\newcommand{\p}{\widetilde{N}}

\begin{document}

\title[Log-Concavity in Narayana Transforms and Stack-Sorting Fibers]
{Log-Concavity in Narayana Transforms \\  and Stack-Sorting Fibers}
\author[]{Colin Defant}
\address[]{Axiom Math, 124 University Avenue, Palo Alto, CA 94301, USA}
\email{colindefant@gmail.com}

\begin{abstract}
Given a permutation $\pi$ of size $n$, let $D_\pi(x)$ denote the descent polynomial of the fiber of $\pi$ under West's stack-sorting map. We prove that $D_\pi(x)$ is ultra-log-concave of order $n-1$ and, therefore, log-concave. More generally, we prove that a multivariate refinement of $D_\pi(x)$ recording ascent tops and descent tops is the basis generating polynomial of a laminar matroid. We deduce a result about simplex integration from Br\"and\'en and Huh's theory and apply it to weighted matchings, obtaining a refinement of the log-concavity of $D_\pi(x)$ that incorporates the peak statistic on the stack-sorting fiber of $\pi$. These results imply the log-concavity of certain sums of products of Narayana polynomials over the lattice points of integral binary nestohedra. A complementary argument shows that replacing the monomials of a homogeneous real-stable polynomial with nonnegative coefficients by products of shifted Narayana polynomials produces a real-rooted polynomial; we discuss several applications. The results in this paper were formalized in Lean by AxiomProver. 
\end{abstract}

\maketitle

\section{Introduction}\label{sec:intro}

\subsection{Stack-Sorting Fibers}
For a nonnegative integer $n$, let $[n]=\{1,\ldots,n\}$. Let $\mathcal W$ denote the set of words over the alphabet of positive integers that have no repeated letters. The symmetric group $S_n$ can be seen as the set of words in $\mathcal W$ that use exactly the letters $1,\ldots,n$. Unless stated otherwise, when referring to a permutation $S_n$, we assume $n\geq 1$. We write $\sigma_i$ for the $i$-th letter in a word $\sigma$.

The \dfn{stack-sorting map} is an operator on $\mathcal W$ whose fibers have surprisingly rich enumerative structure. We denote this operator by $s$. First define $s(\varnothing)=\varnothing$, where $\varnothing$ is the empty word. Given a nonempty word $\sigma\in\mathcal W$, we can write $\sigma=LmR$, where $m$ is the largest letter in $\sigma$; we then define \[
 s(\sigma)=s(L)s(R)m.\]
As an example, we have \[s(416352)=s(41)s(352)\,6=s(1)\,4\,s(3)s(2)\,56=143256.\]
West \cite{West} originally defined this map as a deterministic version of the nondeterministic stack-sorting algorithm studied by Knuth \cite{Knuth}. For $\pi\in\mathcal W$, the number $|s^{-1}(\pi)|$ is called the \dfn{fertility} of $\pi$. This statistic has surprising connections with free probability theory \cite{Troupes}, and it leads to intricate questions about which permutations have specific fertilities and which numbers can arise as fertilities \cite{BousquetMelou,Fertilitopes,Troupes}. See \cite{Thesis} for a thorough discussion of the stack-sorting map.

A \dfn{descent} of a permutation $\sigma\in S_n$ is an index $i\in[n-1]$ such that $\sigma_i>\sigma_{i+1}$. Let $\Des(\sigma)$ be the set of descents of $\sigma$, and let $\des(\sigma)=|\Des(\sigma)|$. For $\pi\in S_n$, we are interested in the polynomial \[D_\pi(x):=\sum_{\sigma\in s^{-1}(\pi)}x^{\des(\sigma)}.\] When $\pi$ is the identity permutation $12\cdots n$, the polynomial $xD_{12\cdots n}(x)$ coincides with the Narayana polynomial
\[N_n(x)=\sum_{j=1}^{n}\frac1n\binom nj\binom n{j-1}x^j.\] For general permutations $\pi$, the polynomials $D_\pi(x)$ are much more complicated. Nevertheless, B\'ona proved that they are always symmetric and unimodal \cite{BonaStack}, and Br\"and\'en established the even stronger fact that these polynomials are $\gamma$-nonnegative \cite{Actions}. The current author subsequently gave an alternative proof of $\gamma$-nonnegativity using combinatorial objects called \emph{valid hook configurations}, and
he conjectured that every polynomial $D_\pi(x)$ is real-rooted \cite{Counting}. He also characterized the polynomials of the form $D_\pi(x)$ via Narayana transforms and the lattice points of certain polytopes called \emph{integral binary nestohedra} (also called \emph{fertilitopes}) \cite{Fertilitopes}.

A nonnegative sequence $(b_j)_{j=0}^d$ is \dfn{log-concave} if $b_j^2\geq b_{j-1}b_{j+1}$ for $1\leq j\leq d-1$, and it has \dfn{no internal zeros} if the indices of its positive terms form an interval. A polynomial with nonnegative coefficients is \dfn{log-concave} if its coefficient sequence is log-concave. A polynomial $\sum_{j=0}^db_jx^j$ is called \dfn{ultra-log-concave of order $d$} if the sequence
$(b_j/\binom{d}{j})_{j=0}^d$ is log-concave and has no internal zeros. Ultra-log-concave polynomials are log-concave.

Until now, it was not even known that the polynomials $D_\pi(x)$ are log-concave. We prove that they are, in fact, ultra-log-concave.

\begin{theorem}\label{thm:fiber}
For each $\pi\in S_n$, the polynomial $D_\pi(x)$ is ultra-log-concave of order $n-1$. In particular, $D_\pi(x)$ is log-concave.
\end{theorem}

For $\sigma\in S_n$ and $i\in[n-1]$, the entry $\sigma_i$ is a \dfn{descent top} of $\sigma$ if $\sigma_i>\sigma_{i+1}$. The entry $\sigma_{i+1}$ is an \dfn{ascent top} of $\sigma$ if $\sigma_i<\sigma_{i+1}$. Let \[
 \DT(\sigma)=\{\sigma_i:\sigma_i>\sigma_{i+1}\}\quad\text{and}\quad  \AT(\sigma)=\{\sigma_{i+1}:\sigma_i<\sigma_{i+1}\}\] denote the sets of descent tops of $\sigma$ and ascent tops of $\sigma$, respectively. Let $\xx=(x_1,x_2,\ldots)$ and $\yy=(y_1,y_2,\ldots)$ be infinite tuples of variables. The key to proving \cref{thm:fiber} is to study the multivariate polynomial
\begin{equation}\label{eq:topPolynomial}
 \cD_\pi(\xx,\yy)= \sum_{\sigma\in s^{-1}(\pi)}\prod_{a\in\DT(\sigma)}x_a\prod_{b\in\AT(\sigma)}y_b.
\end{equation}
Setting $x_i=x$ and $y_i=y$ for all $i$ in $\cD_\pi(\xx,\yy)$ yields the bivariate polynomial
\begin{equation}\label{eq:homogeneousFiber}
 H_\pi(x,y)=\sum_{\sigma\in s^{-1}(\pi)}
 x^{\des(\sigma)}y^{n-1-\des(\sigma)}=y^{n-1}D_\pi(x/y),
\end{equation} which specializes to $D_\pi(x)$ if we set $y=1$.

\begin{theorem}\label{thm:top}
For every permutation $\pi\in S_n$ in the image of $s$, the polynomial $\cD_\pi(\xx,\yy)$ is the basis generating polynomial of a laminar matroid of rank $n-1$. In particular, $\cD_\pi(\xx,\yy)$ is Lorentzian.
\end{theorem}

Lorentzian polynomials, introduced by Br\"and\'en and Huh \cite{BH}, connect coefficient inequalities with discrete convexity and matroid theory. For homogeneous polynomials, the Lorentzian property is equivalent to complete log-concavity, which was developed independently by Anari, Oveis Gharan, and Vinzant \cite{AOGV} and by Anari, Liu, Oveis Gharan, and Vinzant \cite{ALOV}. Every matroid basis generating polynomial is Lorentzian by \cite[Theorem~3.10]{BH}. The contribution of \cref{thm:top} is the identification of the specific matroid behind a stack-sorting fiber. Its independent sets are described by bounds on sizes of subtrees of a plane forest. This description also shows that a permutation in a prescribed stack-sorting fiber is uniquely determined by its set of ascent tops and its set of descent tops.

\begin{figure}[ht]
\includegraphics[height=8.972cm]{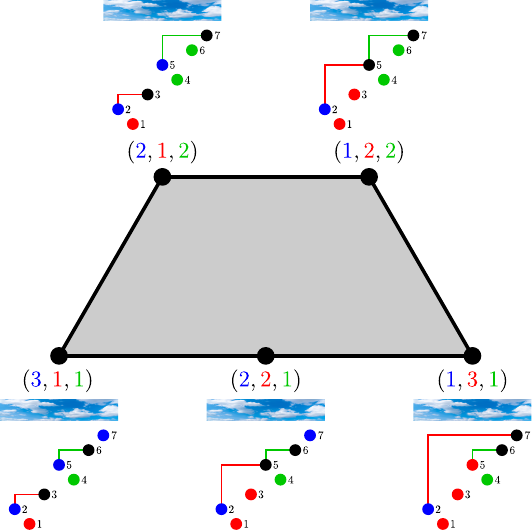}
\caption{The permutation $\pi=2135467$ has $5$ valid hook configurations, which correspond to the $5$ lattice points of the fertilitope $\Fer_\pi$ shown in grey. }
\label{fig:intro}
\end{figure}

\subsection{Narayana Transforms and Weighted Fibers}
Choose $\pi\in S_n$, and let $k=\des(\pi)$. The author introduced a set $\mathrm{VHC}(\pi)$ of combinatorial objects called \emph{valid hook configurations} of $\pi$. Each valid hook configuration of $\pi$ \emph{induces} an integer composition of $n-k$ into $k+1$ parts, called a \dfn{valid composition} of $\pi$. Different valid hook configurations of $\pi$ induce different valid compositions. Let $\cV(\pi)$ be the set of valid compositions of $\pi$. The author proved \cite[Corollary~2.4]{Fertilitopes} that
\begin{equation}\label{eq:introFertility}
 xD_\pi(x)=\sum_{(q_0,\ldots,q_k)\in\cV(\pi)}\prod_{i=0}^{k}N_{q_i}(x).
\end{equation}
For example, the fertilitope of the permutation $\pi=2135467$ is shown in \cref{fig:intro}. Letting $N_{(q_0,\ldots,q_k)}(x)=N_{q_0}(x)\cdots N_{q_k}(x)$, we deduce from \eqref{eq:introFertility} that
\[xD_{2135467}(x)=N_{(3,1,1)}(x)+N_{(2,2,1)}(x)+N_{(1,3,1)}(x)+N_{(2,1,2)}(x)+N_{(1,2,2)}(x)=5x^3+12x^4+5x^5.\]

A multiset $\mathscr A$ of subsets of a set $E$ is \dfn{laminar} if any two sets in $\mathscr A$ that intersect are nested. An \dfn{integral binary nestohedron} is an integral polytope in $\RR_{>0}^r$ obtained as a Minkowski sum of coordinate simplices indexed by a laminar multiset of intervals of $[r]$. The author showed that $\cV(\pi)$ is exactly the set of lattice points of its convex hull $\Fer_\pi$, called the \dfn{fertilitope} of $\pi$. He also showed that the class of nonempty fertilitopes is exactly the same as the class of integral binary nestohedra \cite[Theorems~3.1 and~3.2]{Fertilitopes}. Consequently, \cref{thm:fiber} proves that
\begin{equation}\label{eq:introPolytope}
 \sum_{(q_1,\ldots,q_r)\in P\cap\ZZ^r}\prod_{i=1}^rN_{q_i}(x)
\end{equation}
is ultra-log-concave of order $S+r-1$ for every integral binary nestohedron $P\subseteq\RR_{>0}^r$, where $S$ is the common coordinate sum of points in $P$.

We will give a second proof of \cref{thm:fiber} that accommodates other weights. It uses a polynomial recording unfinished matching edges in a marked forest. Dividing the term corresponding to $h$ unfinished edges by $h!$ turns the forest recurrence into continuous convolution and differentiation. We require several results regarding operations that preserve the Lorentzian property; these follow from results of Br\"and\'en and Huh
\cite[Theorems~3.2 and~3.10]{BH} and results of Ross, S\"uss, and Wannerer \cite{RSW}. Our contribution here is the recurrence for marked forests and its application to weighted stack-sorting fibers and Narayana sums on fertilitopes. The recurrence is controlled by a quadratic kernel whose Lorentzian property has an explicit matrix criterion. The following result illustrates the extra freedom provided by this method.  A \dfn{peak} of $\sigma$ is an integer $i$ with $2\leq i\leq n-1$ such that $\sigma_{i-1}<\sigma_i>\sigma_{i+1}$. Let $\pk(\sigma)$ denote the number of peaks of $\sigma$.

\begin{theorem}\label{thm:peaks}
For every $\pi\in\SSS_n$ and every fixed $\lambda\in[0,2]$, the polynomial
\[H_\pi^\lambda(x,y)= \sum_{\sigma\in s^{-1}(\pi)} \lambda^{\pk(\sigma)} x^{\des(\sigma)}y^{n-1-\des(\sigma)}\] is Lorentzian.
\end{theorem}

\cref{thm:peaks} yields log-concavity results regarding several classes of colored binary plane trees with a fixed postorder. The relevant framework is the theory of \emph{troupes} developed in \cite{Troupes}.

\subsection{A Real-Rooted Narayana Transform}
For each nonnegative integer $m$, let
\[\p_m(x)=N_{m+1}(x)/x.\] Given an infinite tuple $\zz=(z_1,z_2,\ldots)$ of variables and a tuple $\alpha=(\alpha_1,\ldots,\alpha_r)$ of nonnegative integers, let $\zz^\alpha=z_1^{\alpha_1}\cdots z_r^{\alpha_r}$. Define the \dfn{Narayana transform} of a polynomial $F(\zz)=\sum_\alpha c_\alpha\zz^\alpha$ by
\begin{equation}\label{eq:transformIntro}
 \cN(F)(x)=\sum_\alpha c_\alpha\prod_i \p_{\alpha_i}(x).
\end{equation}
A polynomial with real coefficients is \dfn{real-stable} if it does not vanish when every variable has positive imaginary part.

\begin{theorem}\label{thm:stableTransform}
If $F$ is a nonzero homogeneous real-stable polynomial with nonnegative coefficients, then $\cN(F)(x)$ has only real zeros, all of which are negative.
\end{theorem}

The proof of \cref{thm:stableTransform} converts the Narayana polynomials into a signed real-rooted Appell sequence. The algebraic symbol of the resulting operator factors, so a theorem of Borcea and Br\"and\'en \cite{BB} applies.

There are several results in the same spirit as \cref{thm:stableTransform} in the literature. Chen, Wang, and Yang \cite{CWY} proved that the sequence of Narayana polynomials is coefficientwise log-convex; more precisely, for every $m\geq 1$, the polynomial $\p_{m-1}(x)\p_{m+1}(x)-\p_m(x)^2$ has nonnegative coefficients. They also proved that the linear transformation defined by the Narayana number triangle preserves log-convexity of positive sequences. These results concern inequalities between different members of a sequence. Our results concern the zeros and coefficients of individual polynomials obtained by summing products of Narayana polynomials. Mao and Wang proved that the univariate transformation
${x^m\mapsto \p_m(x)}$ preserves real-rootedness on polynomials with
nonnegative coefficients \cite{MaoWang}. The polynomial in \cref{thm:stableTransform} is multivariate, and each monomial is replaced by a product of Narayana polynomials. Real-stable multivariate refinements of Narayana generating polynomials have also been constructed by Yang and Zhang \cite{YangZhang}.

As an example, let $G$ be a graph with vertex set $[r]$ and edge set $E(G)$, and let
$\aa\in\ZZ_{>0}^r$. The polynomial
\[\prod_{v=1}^{r}z_v^{a_v-1}
 \prod_{\{u,v\}\in E(G)}(z_u+z_v)\]
is homogeneous and real-stable, and it has nonnegative coefficients. Thus, \cref{thm:stableTransform} implies real-rootedness of
\[\sum_{\mathcal O}\prod_{v=1}^{r}\p_{a_v-1+\indeg_{\mathcal O}(v)}(x),\]
where $\mathcal O$ ranges over all orientations of $G$ and $\indeg_{\mathcal O}(v)$ is the number of edges oriented toward $v$. If $G$ is a forest, the orientations are in bijection with the lattice points of $\aa+\sum_{\{u,v\}\in E(G)}\Delta_{\{u,v\}}$.
We therefore obtain a real-rootedness theorem for Narayana sums on tree zonotopes with unit edges. These zonotopes allow overlapping edge sets such as $\{1,2\}$ and $\{2,3\}$, so they provide a class beyond fertilitopes.

\subsection{Further Applications}
We apply \cref{thm:fiber}, through \eqref{eq:introPolytope}, to count relations in the noncrossing partition lattice, restricting the block sizes of the upper partition and recording the number of blocks of the lower partition. Using a formula of Kreweras, we express the resulting generating polynomial as a positive constant times a sum of products of Narayana polynomials over the lattice points of a translated simplex. A determinant-twisted Schur polynomial gives another interpretation of the same sums. These applications have classical special cases. We also prove log-concavity of Narayana product sums over the lattice points of matroid base polytopes translated by positive integer vectors.

\subsection{Outline}

\Cref{sec:prelim} supplies the background on trees, valid hook configurations, polytopes, and polynomial inequalities. In \cref{sec:matching}, we construct the
marked forest associated with a permutation and prove \cref{thm:top,thm:fiber}. \Cref{sec:integration} derives the required integration and convolution principles from known preservation results and applies them to the recurrence for marked forests, leading to \cref{thm:peaks} and the general theorem about quadratic kernels. \Cref{sec:troupes} treats colored trees and troupes. In \cref{sec:stable}, we prove \cref{thm:stableTransform} and outline some of its consequences. \Cref{sec:applications} discusses
noncrossing partitions and Schur transforms. \Cref{sec:replacement} proves log-concavity of Narayana sums on translated matroid base polytopes. \Cref{sec:questions} discusses further questions.

\subsection{Formalization}
The results in this article were formalized in Lean by AxiomProver, an internal model at Axiom Math. For explicit details, see
\begin{center}
\url{https://github.com/AxiomMath/LogNarayana}.
\end{center}
This directory includes a formal challenge file containing these theorems, which can be mechanically verified using the Comparator tool in Lean.

\section{Preliminaries}\label{sec:prelim}
Throughout the article, we let $[n]=\{1,\ldots,n\}$ and let $\NN$ denote the
set of nonnegative integers. For integers $a,b$, let $[a,b]_\ZZ=\{j\in\ZZ:a\leq j\leq b\}$. We write $\alpha=(\alpha_1,\ldots,\alpha_r)$ for a tuple in $\NN^r$ and $\zz=(z_1,\ldots,z_r)$ for a tuple of variables, and we use the same subscript convention for other indexed tuples. Write $|\alpha|=\sum_i\alpha_i$, $\alpha!=\prod_i\alpha_i!$, and $\partial_{\zz}^{\alpha}=\prod_i\partial_{z_i}^{\alpha_i}$. Inequalities between tuples are coordinatewise. The \dfn{support} of a polynomial $f(\zz)=\sum_\alpha c_\alpha\zz^\alpha$ is $\{\alpha:c_\alpha\ne0\}$. We also use the notation $[\zz^\alpha]f(\zz)=c_\alpha$. We omit the subscript on $\partial^\alpha$ when the variables are clear, and we write $\Hess f$ for the matrix of second partial derivatives of $f$. For a finite index set $E$ and $I\subseteq E$, let $\zz^I=\prod_{e\in I}z_e$. Let $e_i$ denote the $i$th standard basis vector, let $\one$ denote the all-ones vector, and let $\one_I$ denote the indicator vector of $I$, with the ambient dimension understood from context.

A polynomial is \dfn{multiaffine} if its degree in each variable
is at most 1.

We regard the zero polynomial as Lorentzian. Assertions about the zeros of a polynomial or the M-convexity of its support are understood to apply only when the polynomial is nonzero.

\subsection{Decreasing Binary Trees}
A \dfn{binary plane tree} is a rooted tree such that \begin{itemize}
\item each child of a vertex is designated as either a \dfn{left child} or a \dfn{right child};
\item each vertex has at most one left child and at most one right child.
\end{itemize}
A vertex with two children is
\dfn{binary}, whereas a vertex with only one child is \dfn{unary}. A vertex with no children is a \dfn{leaf}. A binary plane tree is \dfn{full} if every nonleaf vertex is binary, and it is a \dfn{branch} if it has no binary vertices. A labeling of the vertices of a binary plane tree by distinct integers is \dfn{decreasing} if every parent vertex has a label that is larger than those of its children. If $A$ is a finite set of integers, then a \dfn{decreasing binary plane tree on $A$} is a binary plane tree equipped with a decreasing labeling whose set of labels is $A$. When we do not specify the label set $A$, it is assumed to be $[n]$, where $n$ is the number of vertices.

Suppose $\mathcal T$ is a decreasing binary plane tree. Let $\mathcal T_L$ and $\mathcal T_R$ be the left and right subtrees of $\mathcal T$, respectively, each of which is a decreasing binary plane tree. Let $m$ be the label of the root of $\mathcal T$. The \dfn{inorder reading} and the \dfn{postorder reading} of $\mathcal T$, denoted $\mathcal I(\mathcal T)$ and $\mathcal P(\mathcal T)$ respectively, are defined recursively by
\[\mathcal I(\mathcal T)=\mathcal I(\mathcal T_L)\,m\,\mathcal I(\mathcal T_R)\quad\text{and}\quad\mathcal P(\mathcal T)=\mathcal P(\mathcal T_L)\mathcal P(\mathcal T_R)\,m.\]
The inorder and postorder readings of the empty tree are the empty word. The inorder map is a bijection from the set of decreasing binary plane trees on $[n]$ to $\SSS_n$. The postorder map, on the other hand, is not injective. The stack-sorting map is given by
\begin{equation}\label{eq:treeStack}
 s=\mathcal P\circ \mathcal I^{-1}.
\end{equation}
Thus, $\mathcal I^{-1}$ provides a bijection from a stack-sorting fiber $s^{-1}(\pi)$ to the set of decreasing binary plane trees with postorder $\pi$.

In a decreasing binary plane tree $\mathcal T$, a label is a descent top in the inorder reading precisely when the corresponding vertex has a right child. Similarly, a label is an ascent top of $\mathcal I(\mathcal T)$ precisely when its vertex has a left child. Consequently, peaks of the inorder reading of $\mathcal T$ are in bijection with binary vertices of $\mathcal T$.

\begin{figure}[ht]
\includegraphics[height=5.494cm]{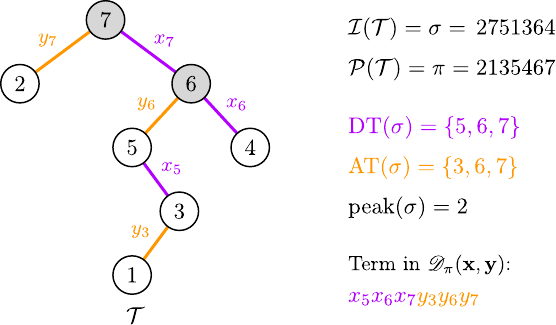}
\caption{A decreasing binary plane tree with inorder reading $\sigma=2751364$ and postorder reading $\pi=2135467$. The two binary vertices in grey correspond to the peaks of $\sigma$.}
\label{fig:treeTops}
\end{figure}

\subsection{Narayana Polynomials}\label{subsec:Narayana}
A polynomial $f(x)=\sum_{j=0}^d b_jx^j$ is \dfn{symmetric} with respect to $d$ if $b_j=b_{d-j}$ for every $j$, and it is \dfn{unimodal} if its coefficients weakly increase and then weakly decrease. It is \dfn{$\gamma$-nonnegative} with respect to $d$ if
\[ f(x)=\sum_{j=0}^{\lfloor d/2\rfloor}\gamma_jx^j(1+x)^{d-2j}, \]
where $\gamma_j\geq 0$ for each $j$.

Let $N_0(x)=1$, and let $N_q(x)=\sum_{j=1}^{q}q^{-1}\binom qj\binom q{j-1}x^j$ for each integer $q\geq 1$. Let $C_j=\frac1{j+1}\binom{2j}{j}$ be the $j$th Catalan number. The shifted Narayana polynomials $\p_m(x)=N_{m+1}(x)/x$ satisfy
\begin{equation}\label{eq:NarayanaGamma}
 \p_m(x)=\sum_{j=0}^{\lfloor m/2\rfloor}
 C_j\binom m{2j}x^j(1+x)^{m-2j}.
\end{equation}
The first few of these polynomials are
\[ \p_0(x)=1,\quad \p_1(x)=1+x,\quad \p_2(x)=1+3x+x^2,\quad \p_3(x)=1+6x+6x^2+x^3. \]
It will be convenient to write
\begin{equation}\label{eq:Qdefinition}
 Q_m(A,B)=\sum_{j=0}^{\lfloor m/2\rfloor} C_j\binom m{2j}B^jA^{m-2j}
\end{equation}
so that $\p_m(x)=Q_m(1+x,x)$. When $A$ is a linear form and $B$ is a quadratic form, $Q_m(A,B)$ is homogeneous of degree $m$.

For every $m\geq 1$, the polynomial $\p_m(x)$ has $m$ distinct negative zeros, which strictly interlace the zeros of $\p_{m+1}(x)$; see \cite[Section~3.4]{LiuWang}.

\subsection{Hooks and Fertilitopes}
The \dfn{plot} of a permutation $\pi=\pi_1\cdots\pi_n$ consists of the points $(i,\pi_i)$. If $i<j$ and $\pi_i<\pi_j$, then the \dfn{hook} with southwest endpoint $(i,\pi_i)$ and northeast endpoint $(j,\pi_j)$ consists of the union of the vertical line segment from $(i,\pi_i)$ to $(i,\pi_j)$ and the horizontal line segment from $(i,\pi_j)$ to $(j,\pi_j)$. As in \cite{Fertilitopes,Troupes}, we define a \dfn{valid hook configuration} of $\pi$ to be a set of hooks on the plot of $\pi$ satisfying the following conditions:
\begin{itemize}
\item the number of hooks is $\des(\pi)$, and the southwest endpoints of the hooks are exactly the points $(i,\pi_i)$ with $i\in\Des(\pi)$;
\item no horizontal segment of a hook lies below a point in the plot of $\pi$;
\item no two hooks intersect except when the northeast endpoint of one is the southwest endpoint of another.
\end{itemize}
Note that the empty set is the unique valid hook configuration of the identity permutation $12\cdots n$.

Consider a valid hook configuration $\mathcal H$ of a permutation $\pi$. Color the hooks in $\mathcal H$ with distinct colors other than blue, and draw a blue sky over the entire plot. This will induce a coloring of the points in the plot of $\pi$ that are not northeast endpoints of hooks. Each of these points receives the color it sees directly above it, looking around the left side of the vertical segment if it is a southwest endpoint of a hook. Let $d_1<\cdots<d_k$ be the descents of $\pi$. The resulting sizes of the color classes form an integer composition
${\bf q}_{\mathcal H}=(q_0,\ldots,q_k)$ of $n-k$. In this tuple, $q_0$ is the number of blue (i.e., sky-colored) points. For each $i\in[k]$, the number $q_i$ counts the points that receive the color of the hook with southwest endpoint $(d_i,\pi_{d_i})$. We call ${\bf q}_{\mathcal H}$ the \dfn{valid composition} of $\pi$ induced by $\mathcal H$. Each valid hook configuration of $\pi$ is determined by the valid composition it induces. Let $\cV(\pi)$ be the set of valid compositions of $\pi$.

\cref{fig:intro} shows the 5 valid hook configurations and valid compositions of the permutation $2135467$, which correspond to the lattice points of the grey trapezoid. Our main use of valid hook configurations and valid compositions comes from \eqref{eq:introFertility}.

For a nonempty subset $J$ of $[r]$, let
$\Delta_J=\conv\{e_j:j\in J\}\subseteq\RR^r$, where $\conv$ denotes convex hull.
A \dfn{binary building set} on $[r]$ is a collection $\cB$ of nonempty intervals containing all singletons such that any two of its members are disjoint or one contains the other. A \dfn{binary nestohedron} is a Minkowski sum
\[P=\sum_{J\in\cB}c_J\Delta_J,\]
where $\cB$ is a binary building set and each coefficient $c_J$ is positive. The polytope $P$ is \dfn{integral} if its vertices have integer coordinates. The singletons in $\cB$ ensure that $P\subseteq\RR_{>0}^r$. The members of $\cB$ can be viewed as selected vertices in a full binary plane tree with leaves $1,\ldots,r$, where we identify a vertex with the set of leaves below it.

Each valid composition of $\pi$ can be seen as a vector in $\mathbb R^{\des(\pi)+1}$. In \cite{Fertilitopes}, the author defined the \dfn{fertilitope} of $\pi$, denoted $\Fer_\pi$, to be the convex hull of $\cV(\pi)$. He proved that $\Fer_\pi\cap\,\ZZ^{\des(\pi)+1}=\cV(\pi)$. (See \cref{fig:intro}.) He also showed that the nonempty fertilitopes are exactly the integral binary nestohedra.

\subsection{Lorentzian Polynomials and Real Stability}\label{subsec:Lorentzian}
A finite nonempty set $J\subseteq\NN^r$ is \dfn{M-convex} if its elements have a common coordinate sum and satisfy the following symmetric exchange property: if $\alpha,\beta\in J$ and $\alpha_i>\beta_i$, then there exists $j$ with $\alpha_j<\beta_j$ such that both $\alpha-e_i+e_j$ and $\beta-e_j+e_i$ belong to $J$. A \dfn{generalized permutohedron} is a polytope all of whose edge directions are parallel to vectors $e_i-e_j$ \cite{Postnikov}.
If $P\subseteq\RR_{\geq 0}^r$ is a nonempty integral generalized permutohedron, then $P\cap\ZZ^r$ is M-convex. Conversely, if $J\subseteq\NN^r$ is a finite nonempty M-convex set, then $\conv(J)$ is an integral generalized permutohedron such that $J=\conv(J)\cap\ZZ^r$.

For our purposes, a homogeneous polynomial $f$ of degree $d\geq2$ with nonnegative coefficients is \dfn{Lorentzian} if its support is M-convex and the Hessian of every derivative $\partial^\alpha f$ with
${|\alpha|=d-2}$ has at most one positive eigenvalue. This is the characterization given in \cite[Theorem~2.25]{BH}.
Nonnegative constants and linear polynomials with nonnegative coefficients are also Lorentzian. For a nonnegative quadratic polynomial, the Hessian condition alone suffices.

A nonnegative function $g$ on a convex set $U\subseteq\RR^r$ is \dfn{log-concave} if
\[g(tu+(1-t)v)\geq g(u)^t g(v)^{1-t}\]
for all $u,v\in U$ and $0<t<1$. When $g$ is positive, this is equivalent to concavity of $\log g$. References to log-concavity on the positive orthant $\RR_{>0}^r$ concern this property of a function; the definition in the introduction concerns coefficient sequences.

We will use the following facts from \cite{BH}. Products, coordinate derivatives, nonnegative linear substitutions, and coefficientwise limits preserve the Lorentzian property. If $f$ is Lorentzian and $\partial^\alpha f$ is nonzero, then $\log(\partial^\alpha f)$ is concave on the positive orthant \cite[Theorem~2.30]{BH}. For a multi-index $\alpha$, write $|\alpha|=\sum_i\alpha_i$
and $\alpha!=\prod_i\alpha_i!$. The normalized support polynomial
\begin{equation}\label{eq:normalizedSupport}
 E_J(\zz)=\sum_{\alpha\in J}\frac{\zz^\alpha}{\alpha!}
\end{equation} is Lorentzian whenever $J$ is M-convex by \cite[Theorem~3.10]{BH}.
Let $\mathbb{I}(\mathsf{M})$ and $\mathbb{B}(\mathsf{M})$ denote the collection of independent sets and the collection of bases, respectively, of a matroid $\mathsf{M}$. All bases have the same cardinality, called the \dfn{rank} of $\mathsf{M}$ and denoted $\rk(\mathsf{M})$. In particular, the basis generating polynomial $B_{\mathsf{M}}(\zz)=\sum_{B\in\mathbb{B}(\mathsf{M})}\prod_{e\in B}z_e$ of any matroid $\mathsf{M}$ is
Lorentzian.

If $\mathscr A$ is a laminar family of subsets of a finite set $E$ and $c\colon\mathscr A\to\NN$ is a function, then
\[\{I\subseteq E:|I\cap A|\leq c(A)\text{ for every }A\in\mathscr A\} \]
is the collection of independent sets of a matroid on $E$ \cite{FifeOxley}. A matroid that admits such a description is called a \dfn{laminar matroid}, and $c(A)$ is called the \dfn{capacity} assigned to $A$.

Let $f(x,y)=\sum_{j=0}^d b_jx^jy^{d-j}$, where $b_j\geq 0$ for all $j$. The Lorentzian property has the following characterization in this case:
\begin{equation}\label{eq:bivariateCriterion}
 \sum_{j=0}^{d}b_jx^jy^{d-j}\text{ is Lorentzian if and only if}\quad \left(\frac{b_j}{\binom dj}\right)_{j=0}^{d} \text{ is log-concave with no internal zeros}.
\end{equation}
For $d\geq2$, the derivatives appearing in the definition are $\partial_x^{j-1}\partial_y^{d-j-1}f$ for $1\leq j\leq d-1$.
Each of these polynomials is homogeneous of degree $2$, and its Hessian has at most one positive eigenvalue if and only if its determinant is nonpositive. These determinant inequalities are exactly the log-concavity inequalities for
$(b_j/\binom dj)_{j=0}^d$. For $f\ne0$, M-convexity of the support is equivalent to the absence of internal zeros. The cases $d=0$ and $d=1$ are immediate.

Let $\halfplane=\{z\in\CC:\operatorname{Im}z>0\}$. A real polynomial $f(z_1,\ldots,z_r)$ is \dfn{real-stable} if $f(z_1,\ldots,z_r)\ne0$ whenever $z_1,\ldots,z_r\in\halfplane$.
A nonzero univariate real polynomial is real-stable if and only if it is real-rooted. Products and identification of variables preserve real stability. Differentiation and specialization of a variable to a real number take real-stable polynomials to polynomials that are either real-stable or zero; see \cite[Lemma~2.4]{Wagner}. By Hurwitz's theorem, if a sequence of real-stable polynomials in a fixed set of variables converges coefficientwise to a nonzero polynomial, then the limit is real-stable; see \cite[Section~2]{Wagner}. A homogeneous real-stable polynomial with nonnegative coefficients is Lorentzian \cite[Proposition~2.2]{BH}. We will use the following result of Borcea and Br\"and\'en.

\begin{theorem}[{\cite[Theorem~1.2]{BB}}]\label{thm:stableSymbol}
Fix $\kappa\in\NN^r$, and let $T$ be a real linear operator on polynomials whose degree in the variable $z_i$ is at most $\kappa_i$ for all $i$. If the polynomial
\begin{equation}\label{eq:algebraicSymbol}
 \Sym_\kappa(T)(\xx,\ww)
 :=\sum_{0\leq\alpha\leq\kappa}\binom{\kappa}{\alpha}
 T(\zz^\alpha)\ww^{\kappa-\alpha}
\end{equation}
is real-stable, then $T$ takes real-stable polynomials to polynomials that are either real-stable or zero.
\end{theorem}

The polynomial in \eqref{eq:algebraicSymbol} is called the \dfn{algebraic symbol} of $T$ for the bounds given by $\kappa$.
Here, $\ww=(w_1,\ldots,w_r)$ is a list of new variables, and $\binom\kappa\alpha=\prod_i\binom{\kappa_i}{\alpha_i}$. In \cref{sec:integration,sec:replacement}, we will also use the following result of Br\"and\'en and Huh \cite[Theorem~3.2]{BH}. Suppose that $T$ is homogeneous of degree
$\ell$, meaning that every nonzero $T(\zz^\alpha)$ is homogeneous of degree $|\alpha|+\ell$. If $\Sym_\kappa(T)(\xx,\ww)$ is Lorentzian, then $T$ sends every homogeneous Lorentzian polynomial
with the prescribed degree bounds to a Lorentzian polynomial.

\section{Marked Forests and Laminar Matroids}\label{sec:matching}

We begin with a matching description of decreasing binary plane trees with
fixed postorder.

\subsection{Matchings}
A \dfn{noncrossing partial matching} on $[N]$ is a set of disjoint pairs $(a,b)$ with ${1\leq a<b\leq N}$ such that no two pairs $(a,b)$ and $(c,d)$ in the set satisfy $a<c<b<d$. We call $a$ and $b$ the left and right endpoints of $(a,b)$.
Vertices that belong to no pair are \dfn{unmatched}. For $\pi\in\SSS_n$, let $\cA(\pi)$ be the set of noncrossing partial matchings on $[n-1]$ such that every element of $\Des(\pi)$ is a left
endpoint and
\begin{equation}\label{eq:weakEdge}
 \pi_a<\pi_{b+1}\qquad\text{for every pair }(a,b).
\end{equation} Define
\[\Gamma_\pi(z)=\sum_{M\in\cA(\pi)}z^{|M|}.\]

\begin{proposition}\label{prop:matching}
For every $\pi\in\SSS_n$, we have
\[D_\pi(x)=
 \sum_{M\in\cA(\pi)} x^{|M|}(1+x)^{n-1-2|M|} =(1+x)^{n-1}\Gamma_\pi\left(\frac{x}{(1+x)^2}\right).\]
More precisely, a matching in $\cA(\pi)$ together with a choice of left or right at each unmatched vertex determines a unique decreasing binary
plane tree with postorder $\pi$.
\end{proposition}

\begin{figure}[htp]
\includegraphics[height=9.124cm]{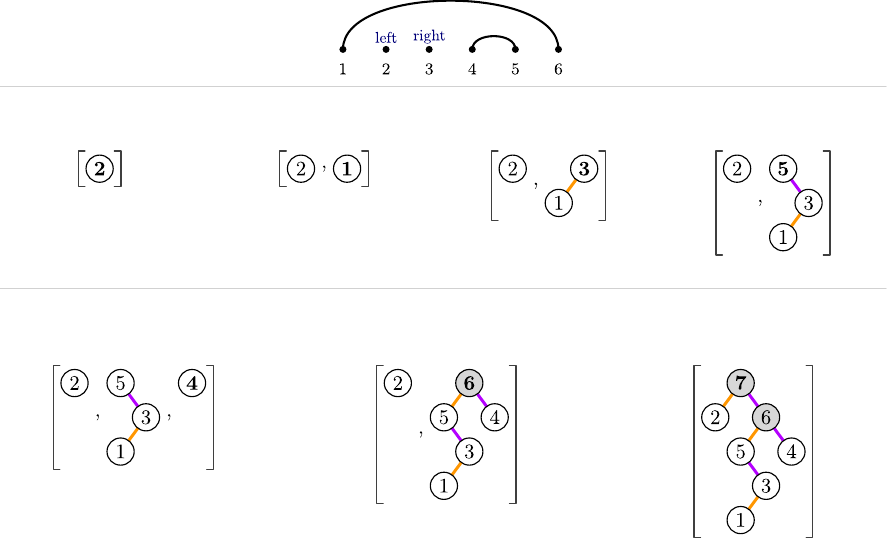}
\caption{The construction in the proof of \cref{prop:matching} for
$\pi=2135467$ and $M=\{(1,6),(4,5)\}$, with the left choice at position $2$ and the right choice at position $3$. The panels show the initial ordered list and the list after each step. The newly introduced label at each step is bold. The final tree is the one in \cref{fig:treeTops}.}
\label{fig:matching}
\end{figure}

\begin{proof}
Fix a matching $M\in\cA(\pi)$, and choose left or right at each unmatched position. We construct a binary plane tree with postorder $\pi$ by successively modifying an ordered list of decreasing binary plane trees (see \cref{fig:matching}). Initially, the list has one entry, which is the 1-vertex tree whose label is $\pi_1$. For $1\leq j\leq n-1$, perform the following operation.
\begin{itemize}
\item If $j$ is a left endpoint of a pair in $M$, append the 1-vertex tree with label $\pi_{j+1}$ to the end of the list.
\item If $j$ is unmatched in $M$, replace the last tree in the list by the tree obtained by adding a new root with label $\pi_{j+1}$.
The old root is the left or right child of the new root according
to the choice at $j$.
\item If $j$ is a right endpoint of $M$, replace the last two trees in the list by a single tree with root label $\pi_{j+1}$, taking the earlier tree as its left subtree and the later tree as its right subtree.
\end{itemize}

To justify the construction, consider $1\leq j\leq n-1$, and let $a_1<\cdots<a_h$ be the left endpoints of the pairs $(a,b)\in M$ satisfying $a<j\leq b$. If $j$ is a right endpoint, the fact that $M$ is noncrossing implies that $(a_h,j)\in M$. Inductively, just before step $j$, the list consists of $h+1$ trees whose root labels, in order, are $\pi_{a_1},\ldots,\pi_{a_h},\pi_j$. Thus, the list is nonempty, and it has at least two entries whenever
the third operation is required.

If $j$ is unmatched or a right endpoint, then $j\notin\Des(\pi)$, so $\pi_j<\pi_{j+1}$. In the third operation, the other child of the new root has label $\pi_{a_h}$, and
$\pi_{a_h}<\pi_{j+1}$ by \eqref{eq:weakEdge}. It follows inductively that every tree in the list is decreasing.
Also, after step $j$, concatenating the postorder readings of the trees in the list gives the word $\pi_1\cdots\pi_{j+1}$.
Each left endpoint increases the length of the list by one, and each right endpoint decreases it by one. Hence, the final list has one entry, which is a decreasing binary plane tree with postorder $\pi$.

Conversely, let $\mathcal T$ be a decreasing binary plane tree with postorder $\pi$. For each binary vertex with label $\pi_{b+1}$, let $\pi_a$ be the label of its left child, and include $(a,b)$ in $M$. Its right child has label $\pi_b$, and its right subtree has postorder $\pi_{a+1}\cdots\pi_b$.
The endpoints of these pairs are distinct because each vertex has at most one parent. The integer intervals $[a+1,b]_{\ZZ}$ belonging to these right subtrees are disjoint or nested, so the pairs form a noncrossing partial matching. Since $\mathcal T$ is decreasing, each
pair satisfies $\pi_a<\pi_{b+1}$.

Unless the vertex labeled $\pi_j$ is the left child of a binary vertex, the next vertex in postorder is its parent and has a larger label. Therefore, every descent position is a left endpoint of $M$, so $M\in\cA(\pi)$. The left and right endpoints are precisely the positions $j$ for which the vertex labeled $\pi_{j+1}$ is a leaf or a binary vertex, respectively. Thus, if $j$ is unmatched, that vertex is unary; choose left or right according to whether its
unique child is a left or right child. The three operations above reconstruct $\mathcal T$ from this matching and these choices, since the subtrees rooted at the children of each vertex appear consecutively in postorder. Conversely, a tree constructed from a matching recovers that matching and its choices by this rule. We have therefore obtained the required bijection.

If $M$ has $m$ pairs, the corresponding trees have $m$ binary vertices and $n-1-2m$ unary vertices. Each binary vertex contributes one descent to the inorder reading. Each unary vertex contributes a
descent exactly when its child is a right child. Thus, by \eqref{eq:treeStack}, the trees corresponding to $M$
contribute $x^m(1+x)^{n-1-2m}$ to $D_\pi(x)$. Summing over $M\in\cA(\pi)$ proves the desired identity.
\end{proof}

\subsection{The Adoptable Forest}
As before, fix $\pi\in S_n$. Call a pair $(a,b)$ with $1\leq a<b\leq n-1$ \dfn{adoptable (with respect to $\pi$)} if
\begin{equation}\label{eq:strongEdge}
 \pi_b=\max\{\pi_{a+1},\ldots,\pi_b\}
 \quad \text{and}\quad
 \pi_{b+1}>\max\{\pi_a,\ldots,\pi_b\}.
\end{equation}
Every pair in a matching in $\cA(\pi)$ is adoptable. Indeed, the entries $\pi_{a+1},\ldots,\pi_b$ are the postorder of the right subtree created by that pair, and $\pi_{b+1}$ is its parent. Conversely, \eqref{eq:strongEdge} implies \eqref{eq:weakEdge}.
We may therefore impose \eqref{eq:strongEdge} in the definition of $\cA(\pi)$ without changing the set.

A \dfn{rooted plane forest} is an ordered list of rooted trees in which the children of each vertex are linearly ordered from left to right. Its \dfn{postorder} visits each vertex after its children, proceeding from left to right through the children and the list of trees.

\begin{lemma}\label{lem:forest}
There is a rooted plane forest $\cF_\pi$ with vertex set $[n-1]$ and postorder reading $12\cdots(n-1)$ such that $(a,b)$ is adoptable with respect to $\pi$ if and only if $b$ is a strict ancestor of $a$ in $\cF_\pi$. Every descent of $\pi$ is a leaf of $\cF_\pi$.
\end{lemma}

\begin{proof}
For fixed $b$, it follows from \eqref{eq:strongEdge} that the left endpoints $a$ for which $(a,b)$ is adoptable form an interval $[\ell_b,b-1]_\ZZ$. If there is no such endpoint, put $\ell_b=b$. Let $I_b=[\ell_b,b]_\ZZ$.

Suppose $(a,c)$ and $(b,d)$ are adoptable and $a<b\leq c<d$. We will show that $(a,d)$ is adoptable. Since $(a,c)$ is adoptable, we have $\pi_c=\max\{\pi_{a+1},\ldots,\pi_c\}$ and $\pi_{c+1}>\max\{\pi_a,\ldots,\pi_c\}$. The adoptability of $(b,d)$ bounds all entries from $b+1$ through $d$ by $\pi_d$. In particular, it bounds $\pi_{c+1}$ by $\pi_d$, so it also bounds $\pi_a,\ldots,\pi_c$ by $\pi_d$. It follows that
$\pi_d=\max\{\pi_{a+1},\ldots,\pi_d\}$ and that $\pi_{d+1}$ exceeds every entry from $\pi_a$ through $\pi_d$. Thus, $(a,d)$ is adoptable.

Now suppose $c\in I_b$ with $c<b$. If $\ell_c<\ell_b$, the fact established in the preceding paragraph can be applied to $(\ell_c,c)$ and $(\ell_b,b)$, contradicting the minimality of $\ell_b$. Therefore, $I_c\subseteq I_b$. This shows that the collection of intervals of the form $I_b$ is laminar. Define the parent of $a$ to be the smallest $b>a$ with $a\in I_b$, if such a $b$ exists. The fact about nesting that we just proved shows that the subtree at $b$ has vertex set exactly $I_b$. The children can be ordered from left to right by their
intervals, giving the required postorder. Finally, if $b\in\Des(\pi)$, no pair ending at $b$ satisfies
\eqref{eq:strongEdge}. It follows that $I_b=\{b\}$ and that $b$ is a leaf.
\end{proof}

We call $\cF_\pi$ the \dfn{adoptable forest} of $\pi$. We mark precisely the leaves in $\Des(\pi)$. Thus, $\cA(\pi)$ consists of the noncrossing matchings on $[n-1]$ in which the right endpoint of each pair is a strict ancestor of its left endpoint in $\cF_\pi$ and every marked leaf is a left endpoint. The remainder of the argument works for any rooted plane forest with any chosen subset of its leaves marked.

\begin{figure}[htp]
\includegraphics[height=4.114cm]{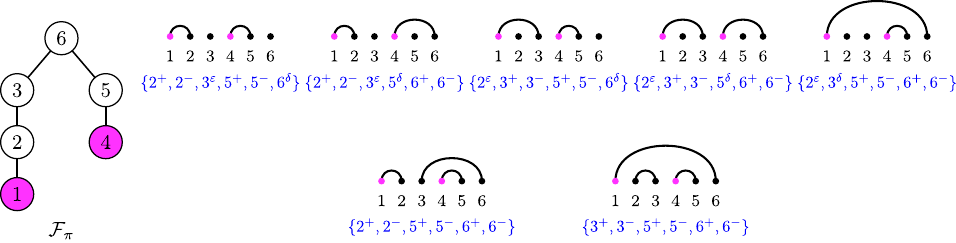}
\caption{The adoptable forest of $\pi=2135467$ along with all 7 matchings in $\cA(\pi)$.
The marked leaves of the forest and the corresponding descent positions $1,4$ in the matchings are pink. Underneath each matching is the corresponding set of bases in the matroid $\mathsf{M}_{\cF_\pi,\Des(\pi)}$. The superscripts $\delta$ and $\varepsilon$ can be chosen arbitrarily from $\{+,-\}$, so each of the matchings in the top row actually corresponds to 4 bases.}
\label{fig:adoptableForest}
\end{figure}

\begin{example}\label{ex:forest}
Let $\pi=2135467$. The forest $\cF_\pi$, shown in
\cref{fig:adoptableForest}, is a tree whose vertex labels are the positions $1,\ldots,6$ and whose marked leaves are $1$ and $4$. For each $b\in[6]$, the subtree rooted at $b$ has vertex set $I_b$, where
\[I_1=\{1\},\quad I_2=[1,2]_\ZZ,\quad I_3=[1,3]_\ZZ,\quad I_4=\{4\},\quad I_5=[4,5]_\ZZ,\quad I_6=[1,6]_\ZZ.\]
The seven matchings in $\cA(\pi)$ are also shown in \cref{fig:adoptableForest}. Five have two pairs, and two have three pairs. Thus, \cref{prop:matching} and \eqref{eq:homogeneousFiber} give
\[ \Gamma_{2135467}(z)=5z^2+2z^3\quad\text{and}\quad H_{2135467}(x,y)=5x^2y^2(x+y)^2+2x^3y^3.
\]
In particular, $D_{2135467}(x)=5x^2+12x^3+5x^4$, in agreement with the calculation from the valid compositions in \cref{fig:intro}.
\end{example}

\subsection{The Laminar Matroid}\label{subsec:laminarMatroid}
Fix a rooted plane forest $\cF$ with $N$ vertices, labeled $1,\ldots,N$ in postorder. Let $L\subseteq[N]$ be a set of leaves of $\cF$, whose elements we call \dfn{marked}. A noncrossing partial matching $M$ on $[N]$ is \dfn{admissible} if for every pair $(a,b)\in M$, the
right endpoint $b$ is a strict ancestor of the left endpoint $a$ in $\cF$, and every vertex in $L$ is a left endpoint of a pair in $M$. For such a matching $M$, let $R(M)$ and $U(M)$ be its set of right endpoints and its set of unmatched vertices, respectively. Define
\begin{equation}\label{eq:forestPolynomial}
 Z_{\cF,L}(\xx,\yy)=  \sum_{M\text{ admissible}}
 \prod_{v\in R(M)}x_vy_v
 \prod_{v\in U(M)}(x_v+y_v).
\end{equation}
The variables are indexed by unmarked vertices. Note that $Z_{\cF,L}(\xx,\yy)$ is homogeneous of degree~$N$.

Let $V_v$ be the vertex set of the subtree of $\cF$ rooted at $v$. Introduce two symbols $v^+,v^-$ for each unmarked vertex $v$, and put
\[
 E=\{v^+,v^-:v\notin L\}\quad\text{and}\quad
 E_v=\{u^+,u^-:u\in V_v\setminus L\}.
\]
Consider the family
\begin{equation}\label{eq:laminarConstraints}
 \mathbb{I}_{\cF,L} =\{I\subseteq E:|I\cap E_v|\leq |V_v|
                 \text{ for every }v\}.
\end{equation}
The sets $E_v$ form a laminar family because the subtree vertex sets $V_v$ are disjoint or nested. If $v\in L$, then $E_v=\varnothing$, so the corresponding inequality in \eqref{eq:laminarConstraints}
is automatic. If $u\notin L$, then $E_u$ contains $u^+$ and $u^-$. It follows that the sets $E_v$ for $v\notin L$ are nonempty and distinct. Thus, assigning capacity $|V_v|$ to each nonempty set $E_v$ realizes $\mathbb{I}_{\cF,L}$ as the collection of independent sets of a laminar matroid. Denote this matroid by $\mathsf{M}_{\cF,L}$. Summing the inequalities in \eqref{eq:laminarConstraints} over the roots of $\cF$ yields the inequality $|I|\leq N$ for every $I\in\mathbb{I}_{\cF,L}$,
so $\rk(\mathsf{M}_{\cF,L})\leq N$.

\begin{theorem}\label{thm:forestMatroid}
If $Z_{\cF,L}(\xx,\yy)\ne0$, then $\mathsf{M}_{\cF,L}$ has rank $N$ and
\[Z_{\cF,L}(\xx,\yy)= \sum_{I\in\mathbb{B}(\mathsf{M}_{\cF,L})} \prod_{v^+\in I}x_v\prod_{v^-\in I}y_v.\]
In particular, $Z_{\cF,L}(\xx,\yy)$ is Lorentzian, and every nonzero coefficient in its monomial expansion is~$1$.
\end{theorem}

\begin{proof}
Fix an admissible matching $M$. To expand its summand in
\eqref{eq:forestPolynomial}, choose $\varepsilon_v\in\{+,-\}$ for each $v\in U(M)$, where $+$ means choosing $x_v$ and $-$ means choosing $y_v$. Associate to these choices the subset
\[I=\{v^+,v^-:v\in R(M)\}\cup\{v^{\varepsilon_v}:v\in U(M)\}
 \subseteq E. \]
All vertices in $R(M)\cup U(M)$ are unmarked, so this subset is
well defined. The resulting monomial is $\prod_{v^+\in I}x_v\prod_{v^-\in I}y_v$, and $|I|=2|M|+|U(M)|=N$.

For a vertex $v$, if $(a,b)\in M$ and $b\in V_v$, then $a\in V_v$ because $a$ is a descendant of $b$ in $\cF$. A pair with both endpoints in $V_v$ contributes two elements to $I\cap E_v$, an unmatched vertex in $V_v$ contributes one, and a left endpoint in $V_v$ whose right
endpoint is outside $V_v$ contributes none. Therefore,
\begin{equation}\label{eq:deficit}
 |V_v|-|I\cap E_v|
 =|\{(a,b)\in M:a\in V_v,\ b\notin V_v\}|.
\end{equation}
Thus, $I\in\mathbb{I}_{\cF,L}$. Since $Z_{\cF,L}(\xx,\yy)\ne0$, at least one admissible matching exists, so this produces an independent set of size $N$. Together with the preceding rank bound, this proves that $\rk(\mathsf{M}_{\cF,L})=N$. Hence, the independent sets of size $N$ are exactly the bases of $\mathsf{M}_{\cF,L}$.

Conversely, fix $I\in\mathbb{I}_{\cF,L}$ with $|I|=N$. Set $c_v=|I\cap\{v^+,v^-\}|$ when $v\notin L$, and set $c_v=0$
when $v\in L$. We construct a matching $M$ and the choices $\varepsilon_v$ by processing the vertices $v$ in postorder. Initially, $M$ and an ordered list of vertices are both empty. The list records vertices that are to be left endpoints and have
not yet been paired. At vertex $v$, perform the following operation.
\begin{itemize}
\item If $c_v=0$, append $v$ to the end of the list.
\item If $c_v=1$, leave $v$ unmatched. Set $\varepsilon_v=+$ if $v^+\in I$, and set $\varepsilon_v=-$ if $v^-\in I$.
\item If $c_v=2$, let $u$ be the last vertex in the list, add $(u,v)$ to $M$, and remove $u$ from the list.
\end{itemize}

We prove inductively that the third operation is always possible
and that the right endpoint of each pair it creates is a strict
ancestor of its left endpoint in $\cF$. Suppose $c_v=2$. The vertices in $V_v\setminus\{v\}$ form a consecutive interval in postorder immediately before $v$. By induction, the right endpoint of each previously formed pair is a strict ancestor of its left endpoint in $\cF$, so a step at a vertex in this interval cannot remove from the list a vertex outside the interval.
Each step at a vertex $u$ changes the length of the list by $1-c_u$. Consequently, just before step $v$, the number of vertices from
$V_v\setminus\{v\}$ still in the list is
\[\sum_{u\in V_v\setminus\{v\}}(1-c_u)
 =|V_v|-1-(|I\cap E_v|-2)=|V_v|+1-|I\cap E_v|\geq 1.\]
All listed vertices outside this interval precede its vertices in
postorder. Thus, the last vertex in the list belongs to
$V_v\setminus\{v\}$, which proves the induction step.

The pairs constructed by this procedure are noncrossing.
Indeed, if $(a,b)$ and $(c,d)$ satisfied $a<c<b<d$, then at step $b$
the vertex $c$ would still be in the list after $a$, so the procedure could not pair $a$ with $b$. At the end, the length of the list is $\sum_{v=1}^N(1-c_v)=N-|I|=0$.
Therefore, every vertex with $c_v=0$ becomes a left endpoint. In particular, every marked leaf becomes a left endpoint, so the resulting matching is admissible.

The matching and choices reconstructed from $I$ select exactly the elements of $I$. Indeed, vertices with $c_v=0,1,2$ become left endpoints, unmatched vertices, and right
endpoints, respectively, and the choices at unmatched vertices record the unique selected elements there. The reconstruction also recovers any admissible matching and choices that produce $I$. At each right endpoint, pairing with any vertex other than the
last vertex in the list would leave a later listed vertex to be paired with a still later right endpoint, creating a crossing.
Hence, the pairing rule is forced, and the two constructions are inverse.

The two inverse constructions prove the claimed identity for $Z_{\cF,L}(\xx,\yy)$ and show that every nonzero coefficient equals $1$. The Lorentzian conclusion follows from the fact that every matroid basis generating polynomial is Lorentzian \cite{BH}.
\end{proof}

\begin{example}
Underneath each matching in \cref{fig:adoptableForest} is the set of bases of $\mathsf{M}_{\cF_\pi,\Des(\pi)}$ associated to it via the correspondence from the proof of \cref{thm:forestMatroid}. The superscripts $\delta$ and $\varepsilon$ can be chosen arbitrarily. For example, the leftmost matching in the top row is labeled by $\{2^+,2^-,3^\varepsilon,5^+,5^-,6^\delta\}$, which represents the 4 bases
\[\{2^+,2^-,3^+,5^+,5^-,6^+\},\, \{2^+,2^-,3^+,5^+,5^-,6^-\},\, \{2^+,2^-,3^-,5^+,5^-,6^+\},\, \{2^+,2^-,3^-,5^+,5^-,6^-\}.\]
\end{example}

\begin{proof}[Proof of \cref{thm:top,thm:fiber}]
Apply \cref{thm:forestMatroid} with $\cF=\cF_\pi$ from
\cref{lem:forest} and $L=\Des(\pi)$. Relabel the variables at position $v$ by replacing $x_v$ with $x_{\pi_{v+1}}$ and replacing $y_v$ with $y_{\pi_{v+1}}$. Under the matching encoding, a left endpoint creates a leaf and contributes $1$; an unmatched vertex creates a unary parent and contributes either $x_{\pi_{v+1}}$ or $y_{\pi_{v+1}}$; a right endpoint creates a binary parent and contributes $x_{\pi_{v+1}}y_{\pi_{v+1}}$. Therefore, the tree interpretation of the descent and ascent top sets identifies \eqref{eq:forestPolynomial} with \eqref{eq:topPolynomial}.
Unused variables can be included as loops in the matroid.
\Cref{thm:forestMatroid} proves \cref{thm:top}.

In $\cD_\pi(\xx,\yy)$, replace each variable $x_a$ with $x$ and each variable $y_a$ with $y$. Nonnegative linear substitution preserves the Lorentzian property, so $H_\pi(x,y)$ is Lorentzian.
Finally, we can apply \eqref{eq:bivariateCriterion} to deduce
\cref{thm:fiber}.
\end{proof}

The polynomial $\cD_\pi(\xx,\yy)$ also yields ultra-log-concavity when specified values are required or forbidden as ascent or descent tops.

\begin{corollary}\label{cor:conditional}
Fix $\pi\in S_n$ and subsets $A_{\mathrm d},B_{\mathrm d},A_{\mathrm a},B_{\mathrm a}\subseteq[n]$. Let $\Omega$ be the set of permutations $\sigma\in s^{-1}(\pi)$
satisfying
\[ A_{\mathrm d}\subseteq\DT(\sigma)\subseteq[n]\setminus B_{\mathrm d} \quad\text{and}\quad A_{\mathrm a}\subseteq\AT(\sigma)\subseteq[n]\setminus B_{\mathrm a}, \]
and suppose $\Omega\ne\varnothing$. Choose positive real numbers $\alpha_a,\beta_a$ for each $a\in[n]$, and give each $\sigma\in\Omega$ the weight
\[\wt(\sigma)=\prod_{a\in\DT(\sigma)}\alpha_a\prod_{b\in\AT(\sigma)}\beta_b.\]
For any subsets $J_{\mathrm d}\subseteq[n]\setminus(A_{\mathrm d}\cup B_{\mathrm d})$ and $J_{\mathrm a}\subseteq[n]\setminus(A_{\mathrm a}\cup B_{\mathrm a})$, the polynomial
\[F(x)=\sum_{\sigma\in\Omega}\wt(\sigma)
 x^{|\DT(\sigma)\cap J_{\mathrm d}|+
     |\AT(\sigma)\cap J_{\mathrm a}|}\] is ultra-log-concave of order $n-1-|A_{\mathrm d}|-|A_{\mathrm a}|$.
\end{corollary}

\begin{proof}
Let $d=n-1-|A_{\mathrm d}|-|A_{\mathrm a}|$. By \cref{thm:top}, the polynomial $\cD_\pi(\xx,\yy)$ is Lorentzian. Let us differentiate respect to each variable $x_a$ for $a\in A_{\mathrm d}$ and each variable $y_b$ for $b\in A_{\mathrm a}$. Then set $x_a=0$ for every $a\in A_{\mathrm d}\cup B_{\mathrm d}$ and $y_b=0$ for every
$b\in A_{\mathrm a}\cup B_{\mathrm a}$. Because $\cD_\pi(\xx,\yy)$ is multiaffine, the resulting polynomial is
\[\sum_{\sigma\in\Omega}
 \prod_{a\in\DT(\sigma)\setminus A_{\mathrm d}}x_a\prod_{b\in\AT(\sigma)\setminus A_{\mathrm a}}y_b.  \]
This polynomial is homogeneous of degree $d$ and Lorentzian.

Now replace each remaining variable $x_a$ by $\alpha_a x$ if $a\in J_{\mathrm d}$ and by $\alpha_a y$ otherwise. Similarly, replace each remaining variable $y_b$ by $\beta_b x$ if $b\in J_{\mathrm a}$ and by $\beta_b y$ otherwise. After multiplying by $\prod_{a\in A_{\mathrm d}}\alpha_a
 \prod_{b\in A_{\mathrm a}}\beta_b$, the resulting polynomial is $y^dF(x/y)$.
These substitutions preserve the Lorentzian property. Thus, $y^dF(x/y)$ is Lorentzian, so \eqref{eq:bivariateCriterion} implies that $F(x)$ is ultra-log-concave of order $d$. \end{proof}

\subsection{Partial Matchings and Factorial Normalization}\label{subsec:partialMatchings}
Fix a rooted plane forest $\cF$ with vertices labeled $1,\ldots,N$ in postorder, and mark a set $L$ of its leaves. The recurrence in \cref{subsec:openRecurrence} will use a generating polynomial that records both a matching and a chosen set of unmatched
vertices. We define this polynomial below and show that dividing the
coefficient of each power $t^h$ by $h!$ produces a Lorentzian polynomial.

Let $M$ be a noncrossing partial matching on $[N]$ such that for each
pair $(a,b)\in M$, the right endpoint $b$ is an ancestor of the
left endpoint $a$ in $\cF$. Let $R(M)$ and $U(M)$ denote the set of right endpoints of $M$ and the set of unmatched vertices of $M$, respectively. Say an unmatched vertex $u$ is \dfn{exposed} if no pair $(a,b)\in M$ satisfies $a<u<b$. Choose a subset $O\subseteq U(M)$ of exposed vertices, which we call \dfn{open vertices}, and require every vertex in $L$ to be either a left endpoint of $M$ or an element of $O$. Define
\[ P_{\cF,L}(\xx,\yy,t)=\sum_{(M,O)}t^{|O|}\prod_{v\in R(M)}x_vy_v \prod_{v\in U(M)\setminus O}(x_v+y_v),\]
where the sum is over all pairs $(M,O)$ satisfying these conditions, and the variables $x_v,y_v$ are indexed by the unmarked vertices. The pairs $(M,\varnothing)$ in this sum are precisely those with $M$ admissible. Thus, $P_{\cF,L}(\xx,\yy,0)=Z_{\cF,L}(\xx,\yy)$.

Recall from \cref{subsec:laminarMatroid} that $E$ consists of two elements $v^+,v^-$ for each unmarked vertex $v$, that $V_v$ is the vertex set of the subtree of $\cF$ rooted at $v$, and that $E_v$ consists of both copies of each unmarked vertex in $V_v$. For a term in the expansion of $P_{\cF,L}(\xx,\yy,t)$, select both $v^+$ and $v^-$ at each right endpoint $v$ of $M$. At each $v\in U(M)\setminus O$, select $v^+$ if $x_v$ is chosen, and select $v^-$ if $y_v$ is chosen. Select neither element at any other unmarked vertex. This defines a subset $I\subseteq E$.
The same count as in \eqref{eq:deficit} gives the identity \[ |V_v|-|I\cap E_v|
 =|\{(a,b)\in M:a\in V_v,\ b\notin V_v\}|+|O\cap V_v|.\] Thus, $I\in\mathbb{I}_{\cF,L}$, the collection of independent sets of $\mathsf{M}_{\cF,L}$, and $|I|=2|M|+|U(M)|-|O|=N-|O|$.

To construct the inverse map, fix $I\in\mathbb{I}_{\cF,L}$, and apply the same argument from the proof of \cref{thm:forestMatroid}, but now without requiring that $|I|=N$. Let $M$ be the resulting matching, and let $O$ be the set of vertices remaining in the list. The same argument proves that every pairing step is possible and that for each resulting pair $(u,v)$, the right endpoint $v$ is a strict ancestor of the left endpoint $u$ in $\cF$. The use of the last listed vertex ensures that $M$ is noncrossing and that every vertex in $O$ is exposed. Note that every marked leaf is either paired as a left endpoint or remains in $O$.

To verify that this reconstruction recovers the original data, consider any pair $(M,O)$ and choices of variables that produce $I$.
At each right endpoint, the pairing rule is forced. Indeed, pairing with
a vertex other than the last listed vertex would leave a later listed vertex that must either be paired still later, creating a crossing, or belong to $O$ and fail to be exposed. Thus, the reconstruction recovers $M$ and then $O$.
At each vertex $v\in U(M)\setminus O$, the choice of $x_v$ or $y_v$ is recovered from which of $v^+$ or $v^-$ belongs to $I$. The two constructions are therefore inverse. Setting $z_{v^+}=x_v$ and $z_{v^-}=y_v$, we obtain
\begin{equation}\label{eq:partialMatroid}
 P_{\cF,L}(\xx,\yy,t)
=\sum_{I\in\mathbb{I}_{\cF,L}}\zz^I t^{N-|I|}.
\end{equation}

For a polynomial $f(\zz,t)=\sum_{h\geq 0}f_h(\zz)t^h$, define its \dfn{factorial normalization in $t$} by
\begin{equation}\label{eq:factorialNormalization}
 \cN_t f(\zz,t)=\sum_{h\geq 0}f_h(\zz)\frac{t^h}{h!}.
\end{equation}
The following proposition is a consequence of the Lorentzian property of matroid basis generating polynomials \cite[Theorem~3.10]{BH}. We include a short proof before applying it to the sum over independent sets for $P_{\cF,L}(\xx,\yy,t)$ in \eqref{eq:partialMatroid}.

\begin{proposition}\label{prop:partialMatroid}
For every matroid $\mathsf{M}$ of rank at most $N$,
\[\sum_{I\in\mathbb{I}(\mathsf{M})}\zz^I\frac{t^{N-|I|}}{(N-|I|)!}\] is Lorentzian.
\end{proposition}

\begin{proof}
Take the direct sum of $\mathsf{M}$ with a free matroid on $k$ new elements, where $k$ is sufficiently large. Then truncate to rank $N$. Specializing each new variable in its basis polynomial to $t/k$ gives the Lorentzian polynomial
\[
 \sum_{I\in\mathbb{I}(\mathsf{M})}\zz^I \binom{k}{N-|I|}\left(\frac tk\right)^{N-|I|}.\]
As $k\to\infty$, this converges coefficientwise to the stated polynomial. Finally, the class of Lorentzian polynomials is closed.
\end{proof}

Since $\mathsf{M}_{\cF,L}$ has rank at most $N$, \cref{prop:partialMatroid} and \eqref{eq:partialMatroid} show that
\[ \cN_t P_{\cF,L}(\xx,\yy,t)
 =\sum_{(M,O)}\frac{t^{|O|}}{|O|!}\prod_{v\in R(M)}x_vy_v\prod_{v\in U(M)\setminus O}(x_v+y_v)\]
is Lorentzian for every rooted plane forest $\cF$ and every set $L$ of marked leaves. Setting $t=0$ recovers $Z_{\cF,L}(\xx,\yy)$. In \cref{sec:integration}, we use simplex integration to extend this conclusion to suitable linear weights on the vertices in $U(M)\setminus O$ and quadratic weights on the pairs in $M$. This extension will prove \cref{thm:peaks}.

\section{Lorentzian Integration and Weighted Forests}\label{sec:integration}

\subsection{Integration over a Simplex}
We first record a theorem that follows from two results of Br\"and\'en and Huh \cite[Theorems~3.2 and~3.10]{BH}. We include the deduction in the integral form needed for the recurrence for marked forests. We use ordinary coordinate measure on a simplex; in particular, for $t\geq 0$, the simplex
\[\{(s_1,\ldots,s_{m-1})\in\RR_{\geq 0}^{m-1}:s_1+\cdots+s_{m-1}\leq t\} \]
has volume $t^{m-1}/(m-1)!$.

\begin{theorem}\label{thm:simplexIntegration}
Let $G(\xx,u_1,\ldots,u_m)$ be a homogeneous Lorentzian polynomial, where $m\geq2$. Then the polynomial
\begin{align}
 \mathfrak I_mG(\xx,t)
 :=\int_{\substack{s_1,\ldots,s_{m-1}\geq 0\\
                  s_1+\cdots+s_{m-1}\leq t}}
 &G(\xx,s_1,\ldots,s_{m-1},  t-s_1-\cdots-s_{m-1})\, \dd s_1\cdots\dd s_{m-1}
 \label{eq:simplexIntegration}
\end{align} is Lorentzian and homogeneous of degree $\deg G+m-1$.
\end{theorem}

\begin{proof}
Write $\xx=(x_1,\ldots,x_r)$ and $\mathbf u=(u_1,\ldots,u_m)$. Integration of monomials gives the identity \begin{equation}\label{eq:DirichletMonomial}
 \mathfrak I_m(\xx^\alpha\mathbf u^\gamma)
 =\frac{\gamma!}{(|\gamma|+m-1)!}\xx^\alpha t^{|\gamma|+m-1}.
\end{equation}
In particular, $\mathfrak I_m$ is a homogeneous linear operator of degree $m-1$.

Choose $\lambda\in\NN^r$ and $\kappa\in\NN^m$ so that the degree of $G$ in $x_j$ is at most $\lambda_j$ for every $j\in[r]$ and its degree in $u_i$ is at most $\kappa_i$ for every $i\in[m]$. Put $k=|\kappa|$, and introduce auxiliary
variables $\yy=(y_1,\ldots,y_r)$ and $\mathbf v=(v_1,\ldots,v_m)$. By \eqref{eq:DirichletMonomial}, the algebraic symbol factors as
\[\Sym_{(\lambda,\kappa)}(\mathfrak I_m)(\xx,t,\yy,\mathbf v)
 =\prod_{j=1}^r(x_j+y_j)^{\lambda_j}\Psi_\kappa(t,\mathbf v), \]
where \[\Psi_\kappa(t,\mathbf v)
=\sum_{0\leq\gamma\leq\kappa}\binom\kappa\gamma\frac{\gamma!}{(|\gamma|+m-1)!}
   t^{|\gamma|+m-1}\mathbf v^{\kappa-\gamma}. \]
Reindexing by $\beta=\kappa-\gamma$ yields \[\Psi_\kappa(t,\mathbf v)=\kappa!\sum_{0\leq\beta\leq\kappa} \frac{t^{k+m-1-|\beta|}\mathbf v^\beta}
        {(k+m-1-|\beta|)!\,\beta!}. \]
Let
\[ J_\kappa=\{(k+m-1-|\beta|,\beta):0\leq\beta\leq\kappa\}.
\] This is the support of the Lorentzian polynomial $t^{m-1}\prod_{i=1}^m(t+v_i)^{\kappa_i}$, so $J_\kappa$ is M-convex. By \cite[Theorem~3.10]{BH}, its normalized support polynomial $E_{J_\kappa}(t,\mathbf v)$ is Lorentzian. Hence, the polynomial $\Psi_\kappa(t,\mathbf v)=\kappa!E_{J_\kappa}(t,\mathbf v)$ is Lorentzian, as is the full symbol.
It follows from \cite[Theorem~3.2]{BH} that $\mathfrak I_mG(\xx,t)$ is Lorentzian. Its degree follows from \eqref{eq:DirichletMonomial}.
\end{proof}

Multiplying the integrand by a monomial in the $u_i$ gives weighted versions of the same result.

For homogeneous polynomials in the tuple of variables $(\xx,t)$, define
\begin{equation}\label{eq:convolution}
 (F*G)(\xx,t)=\int_0^tF(\xx,s)G(\xx,t-s)\dd s.
\end{equation}
Applying \cref{thm:simplexIntegration} to $F(\xx,u_1)G(\xx,u_2)$ gives the result about convolution below.
The result about antidifferentiation was established by Ross, S\"uss, and Wannerer \cite[Corollary~5.9]{RSW}.

\begin{corollary}\label{cor:convolution}
If $F(\xx,t)$ and $G(\xx,t)$ are homogeneous Lorentzian polynomials,
then $(F*G)(\xx,t)$ is Lorentzian.
In particular, the map $F(\xx,t)\mapsto\int_0^tF(\xx,s)\dd s$ preserves the Lorentzian property.
\end{corollary}

The result regarding convolution can also be deduced from \cite[Corollary~3.8]{BH} and \cite[Corollary~5.9]{RSW}, without using \cref{thm:simplexIntegration}. To see this, let $\mathcal E$ denote full factorial normalization, defined on monomials by $\mathcal E(\zz^\alpha)=\zz^\alpha/\alpha!$.
Use disjoint lists $\xx,\yy$, and set
\[Q(\xx,\yy,t)=\mathcal E\bigl(\mathcal E^{-1}(F(\xx,t))\,\mathcal E^{-1}(G(\yy,t))\bigr),\]
where $\mathcal E$ acts on all the variables $(\xx,\yy,t)$. It follows from \cite[Corollary~3.8]{BH} that the polynomial $Q(\xx,\yy,t)$ is Lorentzian. For monomials $\xx^\alpha t^a$ and $\yy^\beta t^b$, the corresponding term in $Q$ is $\frac{a!b!}{(a+b)!}\xx^\alpha\yy^\beta t^{a+b}$.
Consequently,
\[
 \int_0^tQ(\xx,\yy,s)\dd s=\int_0^tF(\xx,s)G(\yy,t-s)\dd s.
\] The left-hand side of the equation is a Lorentzian polynomial by \cite[Corollary~5.9]{RSW}. Identifying $\yy=\xx$ proves the statement about convolution. Thus, both statements of \cref{cor:convolution} follow from existing results.

\subsection{The Open-Endpoint Recurrence}\label{subsec:openRecurrence}
We apply preservation of the Lorentzian property under convolution (\cref{cor:convolution}) and differentiation to weighted versions of
the generating polynomials for partial matchings from \cref{subsec:partialMatchings}. The recurrence below will give \cref{thm:weightedForest} and, through the matching description of stack-sorting fibers, \cref{thm:peaks}.

Let $\cF$ be a rooted plane forest with a subset $L$ of its leaves marked, write $|\cF|$ for its number of vertices, and label these vertices $1,\ldots,|\cF|$ in postorder.
For each pair $(M,O)$ satisfying the conditions in \cref{subsec:partialMatchings}, give each
vertex in $U(M)\setminus O$ weight $A(\xx)$, each pair in $M$ weight $B(\xx)$, and each vertex in $O$ weight $t$, where $A(\xx)$ is homogeneous of degree 1 and $B(\xx)$ is homogeneous of degree 2.
Let $P_{\cF,L}(\xx,t)$ be the sum of the resulting weights over these pairs $(M,O)$. It is homogeneous of degree $|\cF|$.

The polynomial $P_{\cF,L}(\xx,t)$ is the product of the corresponding polynomials for the components of $\cF$. A marked
singleton has polynomial $t$.

Set $P(\xx,t)=P_{\cF,L}(\xx,t)$, and form a rooted plane tree $T$ by adjoining an unmarked root labeled $|\cF|+1$ whose children are the roots of the components of $\cF$, in their given order. Write $P_{\mathrm{root}}(\xx,t)=P_{T,L}(\xx,t)$. For a pair $(M,O)$ on $\cF$, the added root can remain unmatched without being added to $O$, can be added to $O$, or can be paired with a vertex in $O$. In the third case, its partner must be the last vertex of $O$ in postorder and is then removed from $O$. Choosing an earlier vertex of $O$ would enclose a later one within the new pair, contradicting the requirement that open vertices be exposed. Choosing the last one gives a valid matching because that vertex is exposed and is a strict descendant of the added root. Conversely, deleting the added root recovers the pair on $\cF$; if the root was paired, declare its partner open. Therefore,
\begin{equation}\label{eq:ordinaryRoot}
 P_{\mathrm{root}}(\xx,t)=(A(\xx)+t)P(\xx,t)
 +B(\xx)\frac{P(\xx,t)-P(\xx,0)}{t}. \end{equation}
The last term removes one factor of $t$ and inserts the weight of the new pair. There is no factor equal to the number of open vertices because the partner is uniquely determined.

Recall the factorial normalization $\cN_t$ from \eqref{eq:factorialNormalization}. For any polynomials $P(\xx,t)$ and $Q(\xx,t)$, direct integration of monomials gives the identity
\begin{equation}\label{eq:normalizedProduct}
 \cN_t(PQ)=\partial_t\big((\cN_tP)*(\cN_tQ)\big).
\end{equation}
If $U=\cN_tP$, then \eqref{eq:ordinaryRoot} becomes \begin{equation}\label{eq:normalizedRoot}
 \cN_tP_{\mathrm{root}}  =AU+\int_0^tU(\xx,s)\dd s+B\partial_tU.
\end{equation}
The three terms in \eqref{eq:normalizedRoot} should be treated as a single operator.

For a real symmetric matrix $C$, the notation $C\preceq0$ means that $C$ is negative semidefinite. Superscript $\mathsf T$ denotes transpose.

\begin{lemma}\label{lem:kernel}
Let $A(\xx)=\aa^{\mathsf T}\xx$ and $B(\xx)$ have nonnegative
coefficients and degrees 1 and 2, respectively. Set $C=\Hess B$. The kernel
\[ K(\xx,t)=B(\xx)+A(\xx)t+\frac{t^2}{2} \] is Lorentzian if and only if
\begin{equation}\label{eq:kernelCondition}
 C-\aa\aa^{\mathsf T}\preceq0.
\end{equation} When this holds, the operator in \eqref{eq:normalizedRoot} preserves homogeneous Lorentzian polynomials.
\end{lemma}

\begin{proof}
The Hessian of $K(\xx,t)$ is $\left(\begin{smallmatrix}C&\aa\\ \aa^{\mathsf T}&1\end{smallmatrix}\right)$. Completing the square in its quadratic form gives the identity
\[ \xx^{\mathsf T}C\xx+2t\aa^{\mathsf T}\xx+t^2=\xx^{\mathsf T}(C-\aa\aa^{\mathsf T})\xx
  +(t+\aa^{\mathsf T}\xx)^2.\]
The change of variables $(\xx,t)\mapsto(\xx,t+\aa^{\mathsf T}\xx)$ is invertible. Thus, $\Hess K$ has one more positive eigenvalue than $C-\aa\aa^{\mathsf T}$. In particular, $\Hess K$ has at most one positive eigenvalue if and only if $C-\aa\aa^{\mathsf T}\preceq 0$.
The desired equivalence follows from the condition for a quadratic polynomial to be Lorentzian. Moreover,
\begin{equation}\label{eq:kernelIdentity} AU+\int_0^tU(\xx,s)\dd s+B\partial_tU=\partial_t^2(U*K).
\end{equation}
This identity follows by differentiating the convolution twice. Now use \cref{cor:convolution} and closure under differentiation.
\end{proof}

\begin{theorem}\label{thm:weightedForest}
Let $A(\xx)=\aa^{\mathsf T}\xx$ be a linear form with nonnegative coefficients, and let $B(\xx)$ be a homogeneous quadratic polynomial with nonnegative coefficients. Suppose that $\Hess B-\aa\aa^{\mathsf T}$ is negative semidefinite. For every rooted plane forest $\cF$ with a subset $L$ of its leaves marked, the factorially normalized partial-matching polynomial $\cN_tP_{\cF,L}(\xx,t)$ is Lorentzian. Consequently, the generating polynomial for admissible matchings
\[\sum_{M\text{ admissible}}
 B(\xx)^{|M|}A(\xx)^{|\cF|-2|M|}\]
is Lorentzian.
\end{theorem}

\begin{proof}
The empty forest has polynomial $1$, and a marked singleton has polynomial $t$. Joining consecutive forests uses \eqref{eq:normalizedProduct}, which preserves the Lorentzian property by \cref{cor:convolution}. Adjoining an unmarked root uses \eqref{eq:normalizedRoot}, which is covered by \cref{lem:kernel}. Specializing $t=0$ leaves the terms with $O=\varnothing$, which correspond to admissible matchings.
\end{proof}

The proof of \cref{thm:weightedForest} also allows $A$ and $B$ to depend on the unmarked root, provided every corresponding kernel is Lorentzian in the common list of variables. Taking $A_v=x_v+y_v$ and $B_v=x_vy_v$ gives a second proof of the Lorentzian conclusion of
\cref{thm:forestMatroid}.

\begin{figure}[ht]
\includegraphics[height=2.494cm]{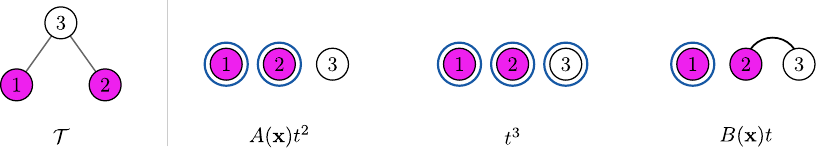}
\caption{A tree $\mathcal T$ (left) and the three terms contributing to the polynomial $P_{\mathrm{root}}({\bf x},t)$ (right).}
\label{fig:root}
\end{figure}

\begin{example}\label{ex:rootRecurrence}
Let $\mathcal F$ consist of two marked singleton trees with labels $1$ and $2$, in this order. The only pair $(M,O)$ on $\mathcal F$ has $M=\varnothing$ and $O=\{1,2\}$, so $P(\mathbf x,t)=t^2$. Adjoin an unmarked root with label $3$ to form the tree $\mathcal T$ in \cref{fig:root}.

The three possible pairs $(M,O)$ on $\mathcal T$ are illustrated to the right of the tree. From left to right, the root $3$ is unmatched and not open, open, or paired with $2$. In the last case, $1$ remains open. Pairing $3$ with $1$ would enclose the open vertex $2$, so it is not allowed. Thus,
\[ P_{\mathrm{root}}(\mathbf x,t)=A(\mathbf x)t^2+t^3+B(\mathbf x)t.\]
Applying factorial normalization yields the identity
\[\mathcal N_tP_{\mathrm{root}}(\mathbf x,t)=A(\mathbf x)\frac{t^2}{2!}+\frac{t^3}{3!}+B(\mathbf x)t.\]
Since $U(\mathbf x,t)=t^2/2!$, the third contribution is $B(\mathbf x)\partial_tU(\mathbf x,t)=B(\mathbf x)t$.
\end{example}

\begin{proof}[Proof of \cref{thm:peaks}]
Take $A=x+y$ and $B=\lambda xy$. The matrix in \eqref{eq:kernelCondition} is
\[\begin{pmatrix}-1&\lambda-1\\ \lambda-1&-1\end{pmatrix},\] whose eigenvalues are $\lambda-2$ and $-\lambda$. Thus, the kernel is Lorentzian exactly when $0\leq\lambda\leq 2$. In the matching encoding, pairs correspond to binary vertices, which are peaks in the inorder reading. Hence, the
generating polynomial in \cref{thm:weightedForest} is exactly $H_\pi^\lambda(x,y)$.
\end{proof}

At $\lambda=1$, this also gives a second proof of
\cref{thm:fiber}. For the running example, \[H_{2135467}^\lambda(x,y)
=5\lambda^2x^2y^2(x+y)^2+2\lambda^3x^3y^3.\] In particular, at $\lambda=2$, the corresponding polynomial after specializing $y=1$ and multiplying by $x$ is
$20x^3+56x^4+20x^5$.

\subsection{Quadratic Replacements on Fertilitopes}
We now apply \cref{thm:weightedForest} to sums over lattice points of integral binary nestohedra. Recall the polynomials $Q_m(A,B)$ from \eqref{eq:Qdefinition}.

\begin{theorem}\label{thm:kernelPolytope}
Let $P\subseteq\RR_{>0}^{r}$ be an integral binary nestohedron. Suppose $A$ and $B$ are homogeneous polynomials with nonnegative coefficients that have degrees 1 and 2, respectively, and satisfy \eqref{eq:kernelCondition}. Then
\begin{equation}\label{eq:kernelPolytope}
 B^{r-1}\sum_{\qq\in P\cap\ZZ^r}\prod_{i=1}^{r}Q_{q_i-1}(A,B)
\end{equation} is Lorentzian. Its degree is $S+r-2$, where $S$ is the common coordinate sum of points in $P$. The polynomial
\[\sum_{\qq\in P\cap\ZZ^r}\prod_{i=1}^{r}\bigl(BQ_{q_i-1}(A,B)\bigr)
\] is also Lorentzian.
\end{theorem}

\begin{proof}
Using the characterization of fertilitopes in \cite[Theorems~3.1 and~3.2]{Fertilitopes}, choose a permutation $\pi$ with $\Fer_\pi=P$. It has $r-1$ descents and has size $S+r-1$. Write $M_m(z)=\sum_j C_j\binom m{2j}z^j$. Combining \eqref{eq:introFertility} with \eqref{eq:NarayanaGamma} and
\cref{prop:matching} yields the identity
\begin{equation}\label{eq:gammaFertility}
 \Gamma_\pi(z)=z^{r-1}\sum_{\qq\in P\cap\ZZ^r}\prod_i M_{q_i-1}(z).
\end{equation}

Weighting the unmatched vertices by $A$ and the pairs by $B$ in \eqref{eq:gammaFertility} gives exactly \eqref{eq:kernelPolytope}. One can verify this first with formal indeterminates $A,B$, so no division by $A$ is needed when $A=0$. Apply \cref{thm:weightedForest}. Finally, \cref{lem:kernel} shows that $B(\xx)+A(\xx)t+t^2/2$ is Lorentzian. Setting $t=0$ shows that $B(\xx)$ is Lorentzian. Since products of Lorentzian polynomials are Lorentzian, multiplying \eqref{eq:kernelPolytope} by $B(\xx)$
proves the last statement.
\end{proof}

\begin{corollary}\label{cor:deformedNarayana}
For $0\leq\lambda\leq2$, put
\[N_q^{(\lambda)}(x) =xQ_{q-1}(1+x,\lambda x)=x\sum_j C_j\binom{q-1}{2j}(\lambda x)^j(1+x)^{q-1-2j}.\]
For every integral binary nestohedron $P\subseteq\RR_{>0}^{r}$, the polynomial
\[\sum_{\qq\in P\cap\ZZ^r}\prod_i N_{q_i}^{(\lambda)}(x)
\] is log-concave and has no internal zeros.
\end{corollary}

\begin{proof}
Let $S$ be the common coordinate sum of the points of $P$. For $\lambda>0$, use $A=x+y$ and $B=\lambda xy$ in \cref{thm:kernelPolytope}, specialize $y=1$, and remove a positive constant and a common monomial factor. These operations preserve ordinary log-concavity and the absence of internal zeros. For $\lambda=0$, the sum is
\[|P\cap\ZZ^r|x^r(1+x)^{S-r}.\]
Its nonzero coefficients form a positive multiple of the binomial sequence $(\binom{S-r}{j})_{j=0}^{S-r}$, which is log-concave.
The polynomial also has no internal zeros.
\end{proof}

Let $e_j(z_1,\ldots,z_r)=\sum_{I\subseteq[r],\,|I|=j}\prod_{i\in I}z_i$ be the $j$th elementary symmetric polynomial, with $e_0(\zz)=1$ and $e_j(\zz)=0$ for $j>r$. A multivariate example comes from setting $A=x_1+\cdots+x_m$ and
$B=e_2(x_1,\ldots,x_m)$. In this case, $C-\aa\aa^{\mathsf T}=-I$, where $I$ is the $m\times m$ identity matrix, so \cref{thm:kernelPolytope} applies.

\subsection{Higher-Degree Kernels}
The convolution operator in \eqref{eq:kernelIdentity} has a higher-degree version. Let $K(\xx,t)$ be homogeneous of degree $m\geq 1$ and Lorentzian, and write it as
\[K(\xx,t)=\sum_{j=0}^{m}W_j(\xx)\frac{t^{m-j}}{(m-j)!},\]
where each nonzero $W_j(\xx)$ is homogeneous of degree $j$. Define $\mathcal R_K(U)=\partial_t^m(U*K)$. The next proposition is a direct consequence of \cref{cor:convolution}.

\begin{proposition}\label{prop:higherKernel}
The operator $\mathcal R_K$ preserves homogeneous Lorentzian polynomials. Every nonzero image has degree 1 greater than its preimage. Explicitly,
\[\mathcal R_K(U)= W_0(\xx)\int_0^tU(\xx,s)\dd s+W_1(\xx)U +\sum_{j=2}^{m}W_j(\xx)\partial_t^{j-1}U.\]
\end{proposition}

For example,
\begin{equation}\label{eq:elementaryKernel}
 K_m(\zz,t)=\sum_{j=0}^{m}e_j(\zz)\frac{t^{m-j}}{(m-j)!}
\end{equation}
is real-stable and, hence, Lorentzian. Indeed, \[K_m(\zz,t)=\lim_{k\to\infty} e_m(z_1,\ldots,z_r,\underbrace{t/k,\ldots,t/k}_{k\text{ copies}}). \]
For every integer $k\geq m$, the polynomial inside the limit is real-stable. The coefficientwise limit is nonzero, so it is also real-stable by Hurwitz's theorem (see \cref{subsec:Lorentzian}).

\section{Postorder Fibers in Troupes}\label{sec:troupes}

\subsection{Troupes}
\begin{figure}[ht]
\includegraphics[height=2.794cm]{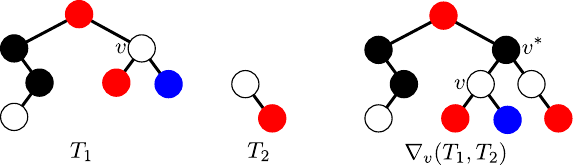}
\caption{Insertion.}
\label{fig:troupe}
\end{figure}

We briefly recall the framework from \cite{Troupes}. Fix a finite set $\mathbf{C}$ of colors, including black and white. A \dfn{colored binary plane tree} is a binary plane tree whose vertices are assigned colors from $\mathbf{C}$. We regard a binary plane tree as a colored binary plane tree in which every vertex is black. Given nonempty colored binary plane trees $T_1,T_2$ and a vertex $v$ of $T_1$, replace $v$ by a left edge with $v$ at the bottom and a new black vertex $v^*$ at the top, and attach $T_2$ as the right subtree of $v^*$. This results in a new colored binary plane tree $\nabla_v(T_1,T_2)$ called the \dfn{insertion} of $T_2$ into $T_1$ at $v$; see \cref{fig:troupe}. A \dfn{troupe} is a family $\mathbf{T}$ of colored binary plane trees satisfying two conditions:
\begin{itemize}
\item Every binary vertex of every tree in $\mathbf{T}$ is black.
\item For all nonempty colored binary plane trees $T_1,T_2$ and every vertex $v$ of $T_1$, we have $T_1,T_2\in\mathbf{T}$ if and only if $\nabla_v(T_1,T_2)\in\mathbf{T}$.
\end{itemize}
See \cite[Definition~2.1]{Troupes}. A statistic $f$ on colored binary plane trees is \dfn{insertion-additive} if we have $f(\nabla_v(T_1,T_2))=f(T_1)+f(T_2)$ for all nonempty colored binary plane trees $T_1,T_2$ and every vertex $v$ of $T_1$.

Given a colored binary plane tree $T$ with $n$ vertices, let $\downarrow \!\!T\!\!\downarrow$ be the set of decreasing colored binary plane trees on $[n]$ that can be obtained by labeling the vertices of $T$. For a troupe $\mathbf{U}$, write $\downarrow\!\!\mathbf{U}\!\!\downarrow=\bigcup_{T\in\mathbf{U}}\downarrow\!\! T\!\!\downarrow$, and let \[D_{\pi,\mathbf{U}}(x)=\sum_{\substack{\mathcal T\in\downarrow \mathbf{U} \downarrow \\ \mathcal P(\mathcal T)=\pi}}x^{\des(\mathcal I(\mathcal T))}.\] When $\mathbf{U}$ is the troupe of binary plane trees, we have $D_{\pi,\mathbf{U}}(x)=D_\pi(x)$. The Refined Tree Fertility Formula \cite[Theorem 4.5]{Troupes} expresses these generating polynomials as product sums over $\cV(\pi)$ whenever the statistics recorded by the weights are insertion-additive. The statistic on colored binary plane trees that records the number of right edges plus $1$ is insertion-additive. On a decreasing labeled tree $\mathcal T$, it equals $\des(\mathcal I(\mathcal T))+1$. We will first prove concrete results directly from \cref{thm:peaks}, then describe the resulting product sums.

Fix positive integers $L,a,b$, and assume $|\mathbf{C}|\geq \max\{L,a,b\}$. Fix sets $X_L, X_a, X_b\subseteq\mathbf{C}$, each containing the color black, such that $|X_j|=j$ for $j\in\{L,a,b\}$. Let $\mathbf{U}_{L,a,b}$ be the set of the colored binary plane trees in which each leaf is assigned a color from the set $X_L$, each unary vertex with a left child is assigned a color from $X_a$, each unary vertex with a right child is assigned a color from $X_b$, and each binary vertex is black. It is straightforward to check that $\mathbf{U}_{L,a,b}$ is a troupe.

\begin{theorem}\label{thm:localColors}
If $L\leq2ab$, then $D_{\pi,\mathbf{U}_{L,a,b}}(x)$ is ultra-log-concave of order $n-1$ for every $\pi\in\SSS_n$ such that $D_{\pi,\mathbf{U}_{L,a,b}}(x)$ is nonzero. In particular, $D_{\pi,\mathbf{U}_{L,a,b}}(x)$ is log-concave.
\end{theorem}

\begin{proof}
Let $\gamma_{\pi,j}=[z^j]\Gamma_\pi(z)$.
A tree with $j$ binary vertices has $j+1$ leaves and $n-1-2j$ unary vertices. Summing over unary orientations and local colors, we find that
\[
\begin{aligned}
 y^{n-1}D_{\pi,\mathbf{U}_{L,a,b}}(x/y)
 &=L\sum_j\gamma_{\pi,j}(Lxy)^j(ay+bx)^{n-1-2j}\\
 &=La^{n-1}H_\pi^{L/(ab)}\left(\frac ba x,y\right).
\end{aligned}
\]
The parameter $L/(ab)$ lies in $[0,2]$. Apply \cref{thm:peaks}, positive rescaling, and \eqref{eq:bivariateCriterion}.
\end{proof}

A \dfn{Schr\"oder two-colored binary tree} is a colored binary plane tree in which every vertex with a left child is black and each other vertex can be black or white. Let $\mathbf{S}$ be the troupe of Schr\"oder two-colored binary trees. For $\pi\in S_n$, define
\[\mathcal C_\pi(z,w)=\sum_{\substack{\mathcal T\in \downarrow\mathbf{S}\downarrow\\\mathcal P(\mathcal T)=\pi}}z^{\black(\mathcal T)}w^{\white(\mathcal T)},\] where $\black(\mathcal T)$ and $\white(\mathcal T)$ denote the numbers of black and white vertices of $\mathcal T$, respectively.

\begin{proposition}\label{prop:Schroeder}
For every $\pi\in\SSS_n$, we have
\[ \mathcal C_\pi(z,w)=(z+w)H_\pi(z+w,z).\] Consequently, $\mathcal C_\pi(z,w)$ is Lorentzian, and the numbers of trees with each possible number of black vertices form an ultra-log-concave sequence of order $n$.
\end{proposition}

\begin{proof}
An underlying decreasing binary plane tree with $d$ right edges has $n-1-d$ vertices with a left child. These vertices are forced to be black. The other $d+1$ vertices can be black or white independently, so the contribution of this tree is $z^{n-1-d}(z+w)^{d+1}$. Summing gives the desired identity. The right side is a nonnegative linear substitution in the Lorentzian polynomial $H_\pi(x,y)$, followed by multiplication by a linear form. Both operations preserve the Lorentzian property.
\end{proof}

For example,
\[\mathcal C_{21345}(z,1)=3z+20z^2+45z^3+42z^4+14z^5.\]
The corresponding factor polynomials are
\begin{equation}\label{eq:SchroederFactors}
 R_q(z)=\sum_{j=0}^{q}\binom{q+j}{q-j}C_jz^j.
\end{equation}
They enumerate the Schr\"oder trees of size $q$ by black vertices \cite[Example 4.12]{Troupes}. Since the map $T\mapsto\black(T)+1$ is insertion-additive, the Refined Tree Fertility Formula tells us that
\[\mathcal C_\pi(z,1)=z^k\sum_{\qq\in\cV(\pi)}\prod_{i=0}^{k}R_{q_i}(z), \]
where $k=\des(\pi)$.  Combining this identity with \cref{prop:Schroeder} and the characterization of fertilitopes in \cite[Theorems~3.1 and~3.2]{Fertilitopes} shows that the sums of products of the polynomials \eqref{eq:SchroederFactors} over integral binary nestohedra are log-concave.

\section{Real-Rooted Narayana Transforms and Tree Zonotopes}\label{sec:stable}

The Lorentzian argument in \cref{sec:integration} proves coefficient inequalities but does not consider real-rootedness. We now use real stability to obtain real-rootedness under a different condition on the indexing polynomial. Recall that $\p_m(x)=N_{m+1}(x)/x$.

\subsection{A Signed Appell Sequence}
An \dfn{Appell sequence} is a sequence $(f_m(t))_{m\geq 0}$ of real polynomials such that $\deg f_m(t)=m$ for every $m\geq 0$ and $f_m'(t)=mf_{m-1}(t)$ for every $m\geq 1$. Define
\begin{equation}\label{eq:signedAppell}
 \eta_m(t)=\sum_{j=0}^{\lfloor m/2\rfloor} (-1)^jC_j\binom m{2j}t^{m-2j}.
\end{equation}
The first few of these polynomials are
\[\eta_0(t)=1,\quad \eta_1(t)=t,\quad \eta_2(t)=t^2-1,\quad\eta_3(t)=t^3-3t,\quad \eta_4(t)=t^4-6t^2+2.\]

\begin{lemma}\label{lem:signedAppell}
The polynomial $\eta_m(t)$ is monic and has $m$ distinct real zeros. Furthermore, we have $\eta_m'(t)=m\eta_{m-1}(t)$ and
\[ \eta_d(z+w)=\sum_{m=0}^{d}\binom dm \eta_m(z)w^{d-m}.\]
\end{lemma}

\begin{proof}
For $r>0$, it follows from \eqref{eq:NarayanaGamma} that
\begin{equation}\label{eq:rootChange}
 \eta_m(r^{-1}-r)=r^{-m}\p_m(-r^2).
\end{equation} The polynomial $\p_m(x)$ has $m$ distinct negative zeros by the classical result recalled in \cref{subsec:Narayana}. Each gives a distinct positive value of $r$ in \eqref{eq:rootChange}. The function $r\mapsto r^{-1}-r$ is a strictly decreasing bijection from $(0,\infty)$ to $\RR$. Hence, $\eta_m(t)$ has $m$ distinct real zeros. Its leading coefficient is $C_0=1$.

Differentiating \eqref{eq:signedAppell} and using the identity $(m-2j)\binom m{2j}=m\binom{m-1}{2j}$ proves that $\eta_m'(t)=m\eta_{m-1}(t)$. Repeated differentiation followed by Taylor's formula proves the final desired identity.
\end{proof}

Consider the real linear operator $\mathfrak A$ on multivariate polynomials defined by
\[\mathfrak A(\zz^\alpha)=\prod_{i=1}^{r}\eta_{\alpha_i}(z_i).\]

\begin{proposition}\label{prop:AppellStability}
The operator $\mathfrak A$ preserves real stability.
\end{proposition}

\begin{proof}
Let $F(\zz)$ be a real-stable polynomial in $z_1,\ldots,z_r$. For each $i\in[r]$, let $\kappa_i$ be the degree of $F(\zz)$ in $z_i$, and put $\kappa=(\kappa_1,\ldots,\kappa_r)$. Restrict $\mathfrak A$ to polynomials whose degree in $z_i$ is at most $\kappa_i$ for every $i\in[r]$. By \eqref{eq:algebraicSymbol}, the algebraic symbol of this restriction is
\begin{align*}
 \Sym_\kappa(\mathfrak A)(\zz,\ww)
 &=\sum_{0\leq\alpha\leq\kappa}\binom\kappa\alpha
   \left(\prod_{i=1}^r \eta_{\alpha_i}(z_i)\right)
   \ww^{\kappa-\alpha}\\
 &=\prod_{i=1}^{r}
   \left(\sum_{m=0}^{\kappa_i} \binom{\kappa_i}{m}\eta_m(z_i)w_i^{\kappa_i-m}\right)\\ &=\prod_{i=1}^{r}\eta_{\kappa_i}(z_i+w_i),
\end{align*}
where the last equality follows from \cref{lem:signedAppell}.
Every factor is real-stable. Indeed, if $z_i,w_i\in\halfplane$, then $z_i+w_i\in\halfplane$, whereas all zeros of $\eta_{\kappa_i}(t)$ are real.
Thus, the product is real-stable, so \cref{thm:stableSymbol} implies that $\mathfrak A(F)$ is real-stable or zero.
Since each $\eta_m(t)$ is monic of degree $m$, the highest homogeneous
component of $\mathfrak A(F)$ equals that of $F$. Thus, $\mathfrak A(F)\ne 0$. Since $F(\zz)$ was arbitrary, this proves the proposition.
\end{proof}

More generally, the proof of \cref{prop:AppellStability} works for any monic real-rooted Appell sequence. Here, the particular coefficients in
\eqref{eq:signedAppell} allow the return to Narayana polynomials.

Let $F(\zz)=\sum_{|\alpha|=D}c_\alpha\zz^\alpha$ be a homogeneous polynomial with nonnegative coefficients, and let
\[\mathcal H_F(t)=\sum_{|\alpha|=D}c_\alpha \prod_i \eta_{\alpha_i}(t),\label{eq:HF}\quad\text{and}\quad
 \Gamma_F(u)=\sum_{|\alpha|=D}c_\alpha
\prod_i M_{\alpha_i}(u),\]
where $M_m(u)=\sum_j C_j\binom m{2j}u^j$. By \eqref{eq:NarayanaGamma} and \eqref{eq:signedAppell}, we have
\[\eta_m(t)=t^mM_m(-t^{-2})\quad\text{and}\quad \p_m(x)=(1+x)^mM_m\left(\frac{x}{(1+x)^2}\right). \]
Since $F(\zz)$ is homogeneous of degree $D$, we have $\sum_i\alpha_i=D$ whenever $c_\alpha\ne0$. Substituting $m=\alpha_i$ in these identities and multiplying over $i\in[r]$ produces the factors $t^D$ and $(1+x)^D$, respectively. Summing with the coefficients $c_\alpha$ gives the identities
\begin{equation}\label{eq:twoChanges}
 \mathcal H_F(t)=t^D\Gamma_F(-t^{-2}) \quad\text{and}\quad \cN(F)(x)=(1+x)^D
                  \Gamma_F\left(\frac{x}{(1+x)^2}\right).
\end{equation}
All these identities are polynomial identities after cancellation of the displayed denominators.

\begin{proof}[Proof of \cref{thm:stableTransform}]
By \cref{prop:AppellStability}, $\mathfrak A(F)$ is real-stable or zero. Identifying its variables gives the real-rooted polynomial $\mathcal H_F(t)$. It is nonzero of degree $D$ since its leading coefficient is $c=\sum_\alpha c_\alpha>0$. Every power of $t$ in \eqref{eq:signedAppell} has exponent $m-2j$, so $\eta_m(-t)=(-1)^m \eta_m(t)$. Since $|\alpha|=D$ in every summand of \eqref{eq:HF}, it follows that $\mathcal H_F(-t)=(-1)^D\mathcal H_F(t)$. Its nonzero zeros consequently occur in opposite pairs. There are (not necessarily distinct) positive numbers $\rho_1,\ldots,\rho_\ell$ such that
\[\mathcal H_F(t)=
 c\,t^{D-2\ell}\prod_{j=1}^{\ell}(t^2-\rho_j^2). \]
The first identity in \eqref{eq:twoChanges} now tells us that $\Gamma_F(u)=c\prod_{j=1}^{\ell}(1+\rho_j^2u)$. Using the second identity yields the explicit factorization
\begin{equation}\label{eq:NarayanaFactorization}
 \cN(F)(x)=c(1+x)^{D-2\ell}
 \prod_{j=1}^{\ell}\big((1+x)^2+\rho_j^2x\big).
\end{equation}
It is straightforward to check that all zeros of this polynomial are negative real numbers.
\end{proof}

The proof of \cref{thm:stableTransform} also shows that every zero of $\Gamma_F(u)$ is negative. The same factorization yields the following weighted extension for all $\lambda\geq 0$.

\begin{corollary}\label{cor:stableDeformation}
If $F$ satisfies the hypotheses of \cref{thm:stableTransform}, then
for every $\lambda\geq 0$, the polynomial
\[ \cN_\lambda(F)(x)=\sum_\alpha c_\alpha
 \prod_i Q_{\alpha_i}(1+x,\lambda x)\]
is real-rooted.
\end{corollary}

\begin{proof}
Replace $\Gamma_F(u)$ by $\Gamma_F(\lambda u)$ in
\eqref{eq:twoChanges}. The polynomial in \eqref{eq:NarayanaFactorization} becomes
\[c(1+x)^{D-2\ell}\prod_j((1+x)^2+\lambda\rho_j^2x),\]
whose zeros are negative for all $\lambda\geq 0$.
\end{proof}

\subsection{Orientation Sums}
Let $G$ be a finite loopless graph with vertex set $[r]$; parallel edges are allowed
and are regarded as distinct. For $\aa\in\ZZ_{>0}^r$, define
\[F_{G,\aa}(\zz)=\prod_{v=1}^{r}z_v^{a_v-1}\prod_{\{u,v\}\in E(G)}(z_u+z_v).
\]
Choosing one endpoint in each edge factor is the same as orienting the edge toward that endpoint. Therefore,
\[F_{G,\aa}(\zz)= \sum_{\mathcal O}\prod_vz_v^{a_v-1+\indeg_{\mathcal O}(v)}. \]
We claim that the polynomial
\begin{equation}\label{eq:orientation}
 \sum_{\mathcal O}\prod_v\p_{a_v-1+\indeg_{\mathcal O}(v)}(x)
\end{equation} is real-rooted. Indeed, each factor of $F_{G,\aa}(\zz)$ is real-stable and has nonnegative coefficients, so we can  apply \cref{thm:stableTransform}.

One can give the two directions of each edge arbitrary nonnegative weights by replacing $z_u+z_v$ by $b_{uv}z_u+c_{uv}z_v$. The same argument gives weighted orientation sums, provided no edge has both weights zero.

\subsection{Lattice Points of Tree Zonotopes}
Given a forest $T$ with vertex set $[r]$ and a vector $\aa\in\ZZ_{>0}^r$, define the translated graphical zonotope
\begin{equation}\label{eq:treeZonotope}
 P_{T,\aa}=\aa+\sum_{\{u,v\}\in E(T)}\Delta_{\{u,v\}}.
\end{equation}

\begin{lemma}\label{lem:treeLattice}
Let $T$ be a forest on $[r]$. The map
\[\mathcal O\longmapsto
 \aa+(\indeg_{\mathcal O}(1),\ldots,\indeg_{\mathcal O}(r)\]
is a bijection from the set of orientations of $T$ to $P_{T,\aa}\cap\ZZ^r$.
\end{lemma}

\begin{proof}
Every orientation chooses an endpoint of every segment in \eqref{eq:treeZonotope}, so its image is a lattice point. For the converse and uniqueness, induct on the number of edges. An isolated vertex contributes only its translation coordinate. Now suppose $T$ has at least one edge, and choose a leaf $v$ with neighbor $u$. In any expression of $\qq\in P_{T,\aa}$ as a sum of segment points, the coordinate $q_v-a_v$ comes only from the edge $\{u,v\}$ and lies in $[0,1]$. If $\qq$ is integral, it equals $0$ or $1$, which forces that segment to contribute $e_u$ or $e_v$, respectively. Subtract the forced contribution, and delete the leaf and its edge. The remaining coordinates give an integral point of the smaller translated forest zonotope. Induction determines the other edge orientations uniquely.
\end{proof}

\begin{theorem}\label{thm:treeZonotope}
For every forest $T$ on $[r]$ and every $\aa\in\ZZ_{>0}^{r}$, the polynomial
\[\sum_{\qq\in P_{T,\aa}\cap\,\ZZ^r}
 \prod_{v=1}^{r}\p_{q_v-1}(x)\] has only negative real zeros. Equivalently,
$\sum_{\qq\in P_{T,\aa}\cap\ZZ^r}\prod_{v=1}^{r}N_{q_v}(x)$ is real-rooted.
\end{theorem}

\begin{proof}
By \cref{lem:treeLattice}, the first polynomial in the statement is the orientation sum in \eqref{eq:orientation}, so it is real-rooted. Its coefficients are nonnegative, and its constant term is $|P_{T,\aa}\cap\ZZ^r|>0$, so all its zeros are negative. Multiplication by $x^r$ gives the second polynomial.
\end{proof}

\begin{example}\label{ex:pathZonotope}
Let $T$ be the path $1-2-3-4$ and let $\aa=\one$. The eight lattice points of $P_{T,\one}$ are
\[
 \begin{gathered}
 (2,2,2,1),\ (2,2,1,2),\ (2,1,3,1),\ (2,1,2,2),\\
 (1,3,2,1),\ (1,3,1,2),\ (1,2,3,1),\ (1,2,2,2).
 \end{gathered}\]
The Narayana sum in \cref{thm:treeZonotope} is
\[8+28x+28x^2+8x^3=4(1+x)(x+2)(2x+1).\]
\end{example}

\subsection{Real-Stable Valid-Composition Polynomials}
Given a permutation $\pi\in\SSS_n$ with $s^{-1}(\pi)\ne\varnothing$, let $k=\des(\pi)$, and define its \dfn{valid-composition polynomial} to be
\[F_\pi(z_0,\ldots,z_k)=
\sum_{\qq\in\cV(\pi)}\prod_{i=0}^{k}z_i^{q_i-1}.\]
This is homogeneous of degree $n-2k-1$. The Refined Fertility Formula states that
\[D_\pi(x)=x^k\cN(F_\pi)(x).\]

\begin{corollary}\label{cor:stableFiber}
If $F_\pi(z_0,\ldots,z_k)$ is real-stable, then $D_\pi(x)$ is real-rooted.
\end{corollary}

For $\pi=1253467$, we have $F_{1253467}(z_0,z_1)=z_0^2z_1(z_0+z_1)$, which is real-stable, so the criterion applies. On the other hand, referring to \cref{fig:intro}, we find that
\[F_{2135467}(z_0,z_1,z_2)
 =z_0^2+z_0z_1+z_1^2+z_2(z_0+z_1). \]
This polynomial is not real-stable since $F_{2135467}(z_0,1,0)=z_0^2+z_0+1$ is not real-rooted. Nevertheless, $D_{2135467}(x)=x^2(5+12x+5x^2)$ has only real zeros.

\section{Restricted Noncrossing Refinements and Schur Transforms}
\label{sec:applications}

\subsection{Restricting the Upper Partition}
A partition of $[n]$ is \dfn{noncrossing} if there do not exist $a<b<c<d$ such that $a,c$ lie in one block and $b,d$ lie in a different block. Let $\NC(n)$ be the set of noncrossing partitions of $[n]$. For $\sigma,\tau\in\NC(n)$, we write $\sigma\leq\tau$ if every block of $\sigma$ is contained in a block of $\tau$. This makes $\NC(n)$ into a lattice whose minimum element $\hat 0$ is the partition into singletons. Write $|\tau|$ for the number of blocks of $\tau$. For integers $n,k,a\geq 1$ with $n\geq ak$, define
\[R_{n,k}^{(a)}(x)=
 \sum_{\substack{\sigma\leq\tau\text{ in }\NC(n)\\
|\tau|=k,\ |B|\geq a\text{ for every }B\in\tau}}
 x^{|\sigma|}. \]
For $0\leq\lambda\leq2$, define
\[R_{n,k}^{(a),\lambda}(x)=
 \sum_{\substack{\tau\in\NC(n),\ |\tau|=k\\
|B|\geq a\text{ for every }B\in\tau}}
 \prod_{B\in\tau}N_{|B|}^{(\lambda)}(x). \]
For a fixed $\tau$, the interval $[\hat 0,\tau]=\{\sigma\in\NC(n):\sigma\leq\tau\}$ is a product of noncrossing partition lattices, one for each block of $\tau$. Since $N_q(x)$ is the generating polynomial for $\NC(q)$ by number of blocks and $N_q^{(1)}(x)=N_q(x)$, we have $R_{n,k}^{(a),1}(x)=R_{n,k}^{(a)}(x)$.

Kreweras proved that for nonnegative integers $m_1,\ldots,m_n$ satisfying $\sum_jm_j=k$ and $\sum_jjm_j=n$, the number of partitions in $\NC(n)$ with $m_j$ blocks of size $j$ for every $j\in[n]$ is
\[\frac{n!}{(n-k+1)!\prod_jm_j!};\] see \cite{Kreweras} or \cite[Proposition~1.4]{Arizmendi}.
There are $k!/\prod_jm_j!$ distinct orderings of these block sizes. Therefore, grouping the terms of $R_{n,k}^{(a),\lambda}(x)$ by the numbers $m_1,\ldots,m_n$ gives the identity \begin{equation}\label{eq:NCsimplex}R_{n,k}^{(a),\lambda}(x)=
 \frac{n!}{k!(n-k+1)!}
 \sum_{\substack{q_1+\cdots+q_k=n\\q_i\geq a}}\prod_{i=1}^kN_{q_i}^{(\lambda)}(x).
\end{equation}

\begin{theorem}\label{thm:restrictedNC}
For all integers $n,k,a\geq 1$ with $n\geq ak$ and all $\lambda\in[0,2]$, the polynomial $R_{n,k}^{(a),\lambda}(x)$ is log-concave and has no internal zeros. In particular, $R_{n,k}^{(a)}(x)$ is log-concave and has no internal zeros.
\end{theorem}

\begin{proof}
The compositions in \eqref{eq:NCsimplex} are precisely the lattice points of the polytope $P=a\one+(n-ak)\Delta_{[k]}$. Since the singletons and $[k]$ form a laminar family of intervals, $P$ is an integral binary nestohedron. Apply \cref{cor:deformedNarayana}, and multiply by the positive scalar in \eqref{eq:NCsimplex}. Setting $\lambda=1$ yields the final statement.
\end{proof}

\begin{example}
For $n=12$, $k=3$, and $a=3$, we have
\begin{align*}
 R_{12,3}^{(3)}(x)=x^3(&220+4290x+30558x^2+105336x^3+192720x^4\\
                      &+192720x^5+105336x^6+30558x^7+4290x^8+220x^9).
\end{align*}
The indexing polytope in \eqref{eq:NCsimplex} is $(3,3,3)+3\Delta_{[3]}$, which has 10 lattice points.
\end{example}

The uniform statement of \cref{thm:restrictedNC} includes cases whose log-concavity can already be obtained by simpler arguments.
For $a=1$, a classical enumerative result \cite{Edelman} yields the identity \[[x^\ell]R_{n,k}^{(1)}(x)=\frac1n\binom n{k-1}\binom n{\ell-k}\binom n{n-\ell}\] for $k\leq\ell\leq n$. The two varying binomial factors are log-concave, so the conclusion follows immediately.

The case $a=2$ has an independent real-rootedness proof using only Lagrange inversion and standard closure properties of real-stable polynomials. We include the argument to make the relationship with earlier methods precise.

\begin{proposition}\label{prop:NCtwo}
For $n\geq2k$, the polynomial $R_{n,k}^{(2)}(x)$ is real-rooted.
\end{proposition}

\begin{proof}
The generating function $U(z)=\sum_{q\geq 1}N_q(x)z^q$ satisfies
\[U=z(U+1)(U+x)\quad \text{and}\quad U-xz=zU(U+1+x).
\]
The composition sum in \eqref{eq:NCsimplex} with $a=2$ and $\lambda=1$ is $[z^n](U-xz)^k$. Put $m=n-k$. Lagrange inversion applied to $U^k(U+1+x)^k$ gives the identity
\[[z^n](U-xz)^k=\frac{k}{m}[u^{n-2k}](u+1+x)^{k-1}(2u+1+x)(u+1)^m(u+x)^m.\]
The product on the right consists of real-stable linear factors in $u,x$. Extracting the coefficient of $u^{n-2k}$ amounts to differentiating $n-2k$ times with respect to $u$, setting $u=0$, and dividing by $(n-2k)!$. These operations produce a polynomial that is real-stable or zero. The composition sum is nonzero because $n\geq2k$, so the resulting univariate real-stable polynomial is real-rooted.
\end{proof}

\subsection{Determinant-Twisted Schur Polynomials}
Let $s_\mu(z_1,\ldots,z_r)$ be the Schur polynomial of a partition $\mu$ with at most $r$ parts. Define
\[\Phi_{\mu,r}(x)=\cN(s_\mu)(x)=\sum_{\alpha}K_{\mu,\alpha}\prod_{i=1}^r \p_{\alpha_i}(x),\]
where the \dfn{Kostka number} $K_{\mu,\alpha}$ is the number of semistandard Young tableaux of shape $\mu$ and content $\alpha$ with entries in $[r]$.

In the partition $(a+d,a^{r-1})$, the notation $a^{r-1}$ means $r-1$ parts equal to $a$.

\begin{corollary}\label{cor:Schur}
For all $a,d\geq 0$ and $r\geq 1$, the polynomial $\Phi_{(a+d,a^{r-1}),r}(x)$ is log-concave. In particular, $\Phi_{\mu,2}(x)$ is log-concave for every partition with at most 2 parts.
\end{corollary}

\begin{proof}
Every semistandard Young tableau of shape $(a+d,a^{r-1})$ with entries in $[r]$ has $a$ columns of height $r$, each filled with $1,\ldots,r$ from top to bottom. The remaining $d$ boxes lie in the first row and can contain any weakly increasing sequence of entries in $[r]$. Consequently, the determinant-twist identity in this case is $s_{(a+d,a^{r-1})}(\zz)=(z_1\cdots z_r)^a h_d(\zz)$, where $h_d(\zz)=\sum_{|\alpha|=d}\zz^\alpha$ is the complete homogeneous symmetric polynomial of degree $d$. It follows that
\[\Phi_{(a+d,a^{r-1}),r}(x) =\sum_{\substack{b_1+\cdots+b_r=d\\b_i\geq 0}} \prod_i \p_{a+b_i}(x).\]
After multiplication by $x^r$, this is the Narayana sum on $(a+1)\one+d\Delta_{[r]}$. The result follows by applying \cref{cor:deformedNarayana} at $\lambda=1$.
\end{proof}

\Cref{cor:Schur} is another interpretation of the product sums in \cref{thm:restrictedNC}. Indeed, if $M=r(a+1)+d$, then
\begin{equation}\label{eq:SchurNC}
 R_{M,r}^{(a+1)}(x)
 =\frac{M!}{r!(M-r+1)!}\,
 x^r\Phi_{(a+d,a^{r-1}),r}(x).
\end{equation}

The one-row case $a=0$ has a classical expression in terms of $P_d^{(r,r)}(t)$, the Jacobi polynomial of degree $d$ with parameters $r,r$ \cite[Sections~18.5 and~18.16]{DLMF}.
The series $V(z)=\sum_{m\geq 0}\p_m(x)z^m$ satisfies $V=(1+zV)(1+xzV)$, so Lagrange inversion yields
\begin{align}
 \Phi_{(d),r}(x)
 &=\frac{r}{r+d}\sum_{j=0}^{d}
       \binom{r+d}{j}\binom{r+d}{d-j}x^j\nonumber\\
 &=\frac{r}{r+d}(1-x)^d
 P_d^{(r,r)}\left(\frac{1+x}{1-x}\right).
 \label{eq:rowJacobi}
\end{align}
The Jacobi zeros lie in $(-1,1)$ \cite{DLMF}, giving negative zeros in \eqref{eq:rowJacobi}. Moreover, \cref{prop:NCtwo} and \eqref{eq:SchurNC} with $a=1$ imply that $\Phi_{(d+1,1^{r-1}),r}(x)$ is real-rooted.

For a general shape $\mu$, some Kostka numbers $K_{\mu,\alpha}$ can exceed $1$. Thus, replacing these coefficients by $1$ and summing only over the monomial support of $s_\mu(\zz)$ would change $\Phi_{\mu,r}(x)$. Huh, Matherne, M\'esz\'aros, and St.~Dizier \cite[Theorem~3]{HMMS} proved that the factorially normalized Schur polynomial $\sum_\alpha K_{\mu,\alpha}\zz^\alpha/\alpha!$ is Lorentzian. Their theorem does not directly establish log-concavity of $\Phi_{\mu,r}(x)$, and our results do not settle this question for arbitrary $\mu$.

\section{Narayana Sums on Matroid Base Polytopes}\label{sec:replacement} We now use the consecutive interlacing of the shifted Narayana polynomials to obtain log-concavity for another class of Narayana sums. The argument is an application of Br\"and\'en and Huh's result on linear operators \cite[Theorem~3.2]{BH}.

The \dfn{base polytope} of a matroid $\mathsf{M}$ on $[r]$ is
\[P_{\mathsf{M}}=\conv\{\one_B:B\in\mathbb{B}(\mathsf{M})\}.\]
Its lattice points are precisely the vectors $\one_B$ with $B\in\mathbb{B}(\mathsf{M})$. Recall that $\mathbb{B}(\mathsf{M})$ is the collection of bases of $\mathsf{M}$ and that
\[
 B_{\mathsf{M}}(\zz)=\sum_{B\in\mathbb{B}(\mathsf{M})}\prod_{i\in B}z_i
\]
is its basis generating polynomial. Write $\widehat\p_m(x,y)=y^m\p_m(x/y)$.

\begin{proposition}\label{prop:interlacingMatroid}
Let $\mathsf{M}$ be a matroid of rank $k$ on $[r]$, and let $\aa\in\NN^r$. The polynomial
\begin{equation}\label{eq:interlacingMatroid}
 \sum_{B\in\mathbb{B}(\mathsf{M})} \prod_{i\in B}\widehat\p_{a_i+1}(x,y) \prod_{i\notin B}\widehat\p_{a_i}(x,y)
\end{equation}
is Lorentzian of degree $|\aa|+k$. Its specialization at $y=1$ is ultra-log-concave of order $|\aa|+k$. If $B_{\mathsf{M}}(\zz)$ is real-stable, this specialization is real-rooted.
\end{proposition}

\begin{proof}
Fix $m\geq 0$. The consecutive interlacing recalled in \cref{subsec:Narayana} implies
that $\p_{m+1}(x)/\p_m(x)$ maps $\halfplane$ into $\halfplane$. Indeed, its partial-fraction expansion has the form
\[\frac{\p_{m+1}(x)}{\p_m(x)}
 =x+c-\sum_{j=1}^m\frac{c_j}{x-\rho_j},\]
where $c\in\RR$, the numbers $\rho_1,\ldots,\rho_m$ are the zeros of $\p_m(x)$, and $c_1,\ldots,c_m>0$. Thus, the polynomial $\p_{m+1}(x)+u\p_m(x)$ is real-stable. Because its coefficients are nonnegative, its homogenization $\widehat\p_{m+1}(x,y)+u\widehat\p_m(x,y)$ is also real-stable; see \cite[proof of Corollary~3.9]{BH}.

Define a linear operator $T$ on multiaffine polynomials by
\[T(\zz^B)=
 \prod_{i\in B}\widehat\p_{a_i+1}(x,y)\prod_{i\notin B}\widehat\p_{a_i}(x,y)\]
for every $B\subseteq[r]$. It is homogeneous of degree $|\aa|$, and its algebraic symbol is
\[\Sym_{\one}(T)(x,y,\ww)=
 \prod_{i=1}^r\bigl(\widehat\p_{a_i+1}(x,y)
                  +w_i\widehat\p_{a_i}(x,y)\bigr).\]
This is homogeneous and real-stable with nonnegative coefficients, so it is Lorentzian. By \cite[Theorem~3.2]{BH}, applying $T$ to the Lorentzian polynomial $B_{\mathsf{M}}(\zz)$ proves that \eqref{eq:interlacingMatroid} is Lorentzian. Its degree is $|\aa|+k$, so \eqref{eq:bivariateCriterion} yields the ultra-log-concavity.

If $B_{\mathsf{M}}(\zz)$ is real-stable, then so is $\zz^{\aa}B_{\mathsf{M}}(\zz)$. Its Narayana transform is the specialization of \eqref{eq:interlacingMatroid} at $y=1$, so the final statement follows from \cref{thm:stableTransform}.
\end{proof}

The lattice points of $\one+\aa+P_{\mathsf{M}}$ are precisely $\one+\aa+\one_B$ for $B\in\mathbb{B}(\mathsf{M})$. Since $N_q(x)=x\p_{q-1}(x)$, the Narayana sum over these points is $x^r$ times the specialization of \eqref{eq:interlacingMatroid} at $y=1$. It is therefore log-concave, and it is real-rooted when
$B_{\mathsf{M}}(\zz)$ is real-stable.

\section{Further Directions}\label{sec:questions}

\subsection{Real-Rootedness of All Fibers}
One of the primary remaining conjectures is the following.

\begin{conjecture}[{\cite{Counting,Fertilitopes}}]\label{conj:fiberRR} For every permutation $\pi$, the nonzero polynomial $D_\pi(x)$ is
real-rooted. Equivalently, the polynomial \[ \sum_{(q_1,\ldots,q_r)\in P\cap\,\ZZ^r}\prod_{i=1}^rN_{q_i}(x)\]
is real-rooted for every integral binary nestohedron $P\subseteq\RR_{>0}^r$.
\end{conjecture}

By \cref{thm:top}, real stability of the basis generating polynomial of a fiber matroid would imply real-rootedness of $D_\pi(x)$. Indeed, setting $x_a=x$ and $y_a=1$ for all $a\in[n]$ in $\cD_\pi(\xx,\yy)$ gives $D_\pi(x)$, and these substitutions preserve real stability. This motivates the following question.

\begin{question}\label{ques:topStable}
For which $\pi\in S_n$ with $s^{-1}(\pi)\ne\varnothing$ is the basis generating polynomial of $\mathsf{M}_{\cF_\pi,\Des(\pi)}$ real-stable?
\end{question}

For $n\geq 1$ and $\pi\in S_n$, the matching description in \cref{prop:matching} gives the identity
\[\sum_{\sigma\in s^{-1}(\pi)}z^{\pk(\sigma)}
 =2^{n-1}\Gamma_\pi(z/4).\]
Thus, the following conjecture asserts log-concavity of the generating polynomial that counts permutations in each fiber $s^{-1}(\pi)$ by their number of peaks.

\begin{conjecture}\label{conj:gamma}
For every $\pi\in S_n$, the polynomial $\Gamma_\pi(z)$ is log-concave.
\end{conjecture}

\cref{conj:gamma} would follow from \cref{conj:fiberRR}. Indeed, by \cref{prop:matching} and the real-rootedness equivalence in \cite[Theorem~1.1]{FPV}, every nonzero $\Gamma_\pi(z)$ would be real-rooted and, hence, log-concave.

\subsection{Beyond Binary Nestohedra}
The author asked whether
\[ \sum_{(q_1,\ldots,q_r)\in P\cap\,\ZZ^r}\prod_{i=1}^rN_{q_i}(x)\]
is real-rooted whenever $P\subseteq\RR_{>0}^r$ is an integral nestohedron, or even whenever it is an integral generalized permutohedron \cite[Question~6.1]{Fertilitopes}.
The corresponding log-concavity problem is also natural. \Cref{thm:treeZonotope} proves real-rootedness for the translated unit-edge forest zonotopes $P_{T,\aa}$, which form a class of integral generalized permutohedra.

\begin{question}\label{ques:dilatedTrees}
Fix a tree $T$ on $[r]$, a tuple $\aa\in\ZZ_{>0}^{r}$, and a positive integer $m_e$ for each edge $e\in E(T)$. Set
\[P=\aa+\sum_{e\in E(T)}m_e\Delta_e\subseteq\RR_{>0}^r. \]
Is the polynomial
\[ \sum_{(q_1,\ldots,q_r)\in P\cap\,\ZZ^r}\prod_{i=1}^rN_{q_i}(x)\]
real-rooted? Is it log-concave?
\end{question}

It turns out that some condition on the polytope is necessary. The following example answers the final part of \cite[Question~6.2]{Fertilitopes}.

\begin{proposition}\label{prop:arbitraryPolytope} There is an integral polytope $P\subseteq\RR_{>0}^{27}$ contained in $\{\qq\in\RR^{27}:\sum_{i=1}^{27}q_i=51\}$ such that the polynomial
\[
 \sum_{\qq\in P\cap\,\ZZ^{27}}\prod_{i=1}^{27}N_{q_i}(x)
\] is not log-concave.
\end{proposition}

\begin{proof}
Set $\aa=(23,1,\ldots,1)\in\ZZ^{22}$ and $\mathbf b=(2,\ldots,2)\in\ZZ^{22}$. For every 2-element subset $S\subseteq[5]$, let $u_S=\one+\one_S\in\ZZ^5$ and define
\[
 v_S=
 \begin{cases}
  (\aa,u_S),&S=\{1,2\};\\
  (\mathbf b,u_S),&S\ne\{1,2\}.
 \end{cases}\]
Let $P=\conv\{v_S:S\subseteq[5],\ |S|=2\}$. All 10 points $v_S$ have positive integer coordinates and coordinate sum $51$. We claim that they are the only lattice points of $P$. Let $\qq\in P\cap\ZZ^{27}$. Its last five coordinates are integers in $\{1,2\}$ with sum $7$, so they equal $\one+\one_{S^*}$ for a unique 2-element subset $S^*\subseteq[5]$. Write $\qq=\sum_{S\subseteq[5],\,|S|=2}\lambda_Sv_S$, where $\lambda_S\geq 0$ and $\sum_S\lambda_S=1$. Comparing the last five coordinates gives $\one_{S^*}=\sum_S\lambda_S\one_S$. For each $i\notin S^*$, the $i$th coordinate on the right must be zero. Hence, $\lambda_S>0$ implies $S\subseteq S^*$, so $S=S^*$. Thus, $\qq=v_{S^*}$, proving the claim. The same argument shows that each $v_S$ is a vertex of $P$.

Consequently, \[x^{-27}\sum_{\qq\in P\cap\,\ZZ^{27}}\prod_{i=1}^{27}N_{q_i}(x)=(1+x)^2\big(\p_{22}(x)+9(1+x)^{22}\big).\]
Its first three coefficients are $10,471,22472$, which violate log-concavity.
\end{proof}

The polytope $P$ constructed in the proof of \cref{prop:arbitraryPolytope} is not a generalized permutohedron. The 9 vertices with first coordinate $2$ form a base face, and the remaining vertex is the apex of a pyramid. An edge from the apex to the base changes all of the first 22 coordinates, so its direction is not parallel to any vector of the form $e_i-e_j$. Thus, \cref{prop:arbitraryPolytope} does not answer the question of whether Narayana product sums over integral generalized permutohedra are always real-rooted.

\section*{AI Usage Statement}
Back in 2021, the author suspected that there might be some way to construct a matroid from a stack-sorting fiber $s^{-1}(\pi)$ in order to find a Lorentzian refinement of the polynomial $D_\pi(x)$. This was the starting point for the current work. The author discussed this problem with AxiomProver, an internal model at Axiom Math. AxiomProver found the matroid $\mathsf{M}_{\cF_\pi,\Des(\pi)}$, which was the main tool needed to prove that $D_\pi(x)$ is log-concave. In a series of further prompts, the author steered AxiomProver toward the other results discussed in this paper (and away from several less interesting results). Once all proofs were in place, the author asked AxiomProver to write the initial draft of the paper. The author then rewrote the entire paper, crafting all figures from scratch and checking all details along the way. The author accepts full responsibility for all results presented in the article. After the article was written, the author used AxiomProver to formalize all of the results in Lean relative to standard existing literature.

\end{document}